\documentclass{amsart}
\usepackage{amsmath}
\usepackage{amsfonts}
\usepackage{amssymb}
\usepackage{graphicx}
\usepackage{amssymb}
\usepackage{amscd}
\usepackage{mathrsfs}
\usepackage{euscript}

\tolerance=1600 \hbadness=10000 \theoremstyle{plain}
\theoremstyle{plain}
\newtheorem{thm}{Theorem}[section]
\newtheorem{prop}[thm]{Proposition}
\newtheorem{cor}[thm]{Corollary}
\newtheorem{lem}[thm]{Lemma}

\theoremstyle{definition}
\newtheorem{notation}[thm]{Notation}
\newtheorem{defn}[thm]{Definition}
\theoremstyle{remark}
\newtheorem{remark}[thm]{Remark}

\newtheorem{example}[thm]{Example}
\renewenvironment{proof}[1][Proof]{\textbf{#1.} }{\hfill \rule{0.5em}{0.5em}}
\numberwithin{equation}{section}

\title{Heat kernel analysis on semi-infinite Lie groups}
\author{Tai Melcher}
\address{Department of Mathematics, University
of Virginia, Charlottesville, VA 22936}
\email{melcher@virginia.edu}


\begin{document}
\begin{abstract}
This paper studies Brownian motion and heat kernel measure on a class of
infinite dimensional Lie groups.  We prove a Cameron-Martin type
quasi-invariance theorem for the heat kernel measure and give estimates on the
$L^p$ norms of the Radon-Nikodym derivatives.   We also prove that
a logarithmic Sobolev inequality holds in this setting.
\end{abstract}

\maketitle
\tableofcontents

\section{Introduction}
\label{s.intro}
We define Brownian motion on a class of infinite dimensional Lie
algebras which we call {\em semi-infinite Lie algebras.}  We then prove a
Cameron-Martin type quasi-invariance result for the associated heat
kernel measure, as well as a logarithmic Sobolev inequality.  A
particular example of these semi-infinite Lie algebras was treated
in \cite{DG08-2}, and we build on the methods used there.

We briefly describe here the main results and give an outline of the paper;
see Sections \ref{s.prelim} and \ref{s.semi-inf} for definitions.
Let $(W,H,\mu)$ be an abstract Wiener space and $\mathfrak{v}$ be
a finite dimensional Lie algebra equipped with an inner product.
Let $\mathfrak{g}=W\oplus\mathfrak{v}$ be a Lie algebra extension
of $W$ by $\mathfrak{v}$, and we will call
$\mathfrak{g}_{CM}=H\oplus\mathfrak{v}$ the Cameron-Martin Lie
subalgebra of $\mathfrak{g}$.  If $\mathfrak{g}$ is nilpotent, we
may define an explicit group operation on $\mathfrak{g}$ via the
Baker-Campbell-Hausdorff-Dynkin formula, and $W\oplus\mathfrak{v}$ 
equipped with
this group operation will be denoted by $G$.  Similarly,
$G_{CM}=H\oplus\mathfrak{v}$ with the same group operation
is called the Cameron-Martin subgroup of $G$, and we equip $G_{CM}$ with the
left invariant Riemannian metric which agrees with the inner product
\[ \langle (A,a),(B,b) \rangle_{\mathfrak{g}_{CM}}
	= \langle A,B \rangle_H + \langle a,b \rangle_\mathfrak{v}
\]
on $\mathfrak{g}_{CM}\cong T_\mathbf{e} G_{CM}$.

In Section \ref{s.prelim}, we set the notation and give some standard
facts needed about abstract Wiener spaces and extensions of Lie
algebras.  In Section \ref{s.semi-inf}, we construct the {\em
semi-infinite Lie algebras} and give some examples.  We make some
additional requirements so that the Lie bracket on $\mathfrak{g}$
is continuous, making $\mathfrak{g}$ into a Banach Lie algebra.  In
Section \ref{s.norm}, this gives bounded Hilbert-Schmidt norms for
the Lie bracket, and, in Section \ref{s.Ric}, lower bounds
on the Ricci curvature of $G$ and a uniform lower bound on certain
finite dimensional approximations of $G$.

In Section \ref{s.BM}, we define Brownian motion on $G$ as the solution to a 
stochastic differential equation with respect to a Wiener process on
$\mathfrak{g}$.  Let $B_t$ denote Brownian motion on $\mathfrak{g}$.
Then, Brownian motion on $G$ is the solution to the Stratonovich
stochastic differential equation
\[ \delta g_t = g_t \,\delta B_t := L_{g_t*} \delta B_t, 
	\quad \text{ with } g_0=\mathbf{e}=(0,0). \]
For $t>0$, let $\Delta_n(t)$ denote the simplex in
$\mathbb{R}^n$ given by
\[ \{s=(s_1,\cdots,s_n)\in\mathbb{R}^n: 0<s_1<s_2<\cdots<s_n<t\}. \]
Let $\mathcal{S}_n$ denote the permutation group on $(1,\cdots,n)$,
and, for each $\sigma\in \mathcal{S}_n$, let $e(\sigma)$ denote the
number of ``errors'' in the ordering
$(\sigma(1),\sigma(2),\cdots,\sigma(n))$, that is, $e(\sigma)=\#
\{j<n: \sigma(j)>\sigma(j+1)\}$.
Then the Brownian motion on $G$ may be written as 
\[
g_t = \sum_{n=1}^{r-1} \sum_{\sigma\in\mathcal{S}_n}  
	\left( (-1)^{e(\sigma)}\bigg/ n^2 
	\begin{bmatrix} n-1 \\ e(\sigma) \end{bmatrix}\right) 
	\int_{\Delta_n(t)} 
	[ [\cdots[\delta B_{s_{\sigma(1)}},\delta B_{s_{\sigma(2)}}],\cdots], 
	\delta B_{s_{\sigma(n)}}],
\]
where this sum is finite since $\mathfrak{g}$ is assumed to be nilpotent.
In Section \ref{s.BM}, we show that these stochastic integrals are
well-defined and each may be expressed as a sum of iterated It\^o
integrals.  We also show that $g_t$ may be realized as a
limit of Brownian motions living on the finite dimensional
approximations to $G$.  In particular, we show in Proposition
\ref{p.approx} that this convergence holds in $L^p$, for all
$p\in[1,\infty)$.

In Theorem \ref{t.quasi}, we apply the previous results and a theorem from
\cite{DG08-1} to prove that $\nu_t=\mathrm{Law}(g_t)$ is invariant
under (right or left) translation by elements of $G_{CM}$.  Moreover,
this theorem gives good bounds on the $L^p$-norms of the Radon-Nikodym
derivatives.  These results are important for future 
applications to spaces of holomorphic functions on $G$, as in \cite{DG08-3}.  
We also show in Theorem \ref{t.logsob}
that a logarithmic Sobolev inequality holds for polynomial cylinder
functions on $G$.

For heat kernel analysis, quasi-invariance results, and 
logarithmic Sobolev inequalities 
in related infinite dimensional settings, see \cite{Malliavin06, Inahama04}.

\section{Preliminaries}
\label{s.prelim}
\subsection{Abstract Wiener spaces}
\label{s.wiener}
In this section, we summarize several well known properties of
Gaussian measures and 
abstract Wiener spaces that are required for the sequel.  For proofs
of these results, see Section 2 of \cite{DG08-2}.  Also see
\cite{Bog98,Kuo75} for more on abstract Wiener spaces and some particular
examples.  

Suppose that $W$ is a real separable Banach space and $\mathcal{B}_{W}$ is
the Borel $\sigma$-algebra on $W$.

\begin{defn}
\label{d.2.1} 
A measure $\mu$ on $(W,\mathcal{B}_{W})$ is called a (mean zero,
non-degenerate) {\it Gaussian measure} provided that its characteristic
functional is given by
\begin{equation}
\label{e.2.1}
\hat{\mu}(u) := \int_W e^{iu(x)} d\mu(x)
	= e^{-\frac{1}{2}q(u,u)}, \qquad \text{ for all } u\in W^*,
\end{equation}
for $q=q_\mu:W^*\times W^*\rightarrow\mathbb{R}$ a symmetric, positive
definite quadratic form.
That is, $q$ is a real inner product on $W^*$.
\end{defn}

\begin{thm}
\label{t.2.3}
Let $\mu$ be a Gaussian measure on a real separable Banach space $W$.  
For $1\le p<\infty$, let
\begin{equation}
\label{e.2.2}
C_p :=\int_W \|w\| _{W}^{p} \,d\mu(w). 
\end{equation}
For $w\in W$, let
\[ 
\|w\|_H := \sup\limits_{u\in W^*\setminus\{0\}}\frac{|u(w)|}{\sqrt{q(u,u)}}
\]
and define the {\em Cameron-Martin subspace} $H\subset W$ by
\[ H := \{h\in W : \|h\|_H < \infty\}. \]
Then
\begin{enumerate}
\item \label{i.1}
For all  $1\le p<\infty$, $C_p<\infty$.

\item $H$ is a dense subspace of $W$.

\item There exists a unique inner product $\langle\cdot,\cdot\rangle_H$ 
on $H$ such that $\|h\|_H^2 = \langle h,h\rangle_H$ for all $h\in H$, and 
$H$ is a separable Hilbert space with respect to this inner product.

\item \label{i.3}
For any $h\in H$,
$\|h\|_W \le \sqrt{C_2} \|h\|_H$.

\item \label{i.5}
If $\{k_j\}_{j=1}^\infty$ is an orthonormal basis of $H$ 
and $\varphi$ is a bounded linear map from $W$ to a real Hilbert space
$\mathbf{C}$, then
\begin{equation}
\label{e.2.13} 
\| \varphi\| _{H^*\otimes\mathbf{C}}^{2}
	:=\sum_{j=1}^\infty \|\varphi(k_j)\|_\mathbf{C}^2
	= \int_W \|\varphi(w)\|_\mathbf{C}^2 \,d\mu(w)  
	< \infty. 
\end{equation}
\end{enumerate}
\end{thm}

A simple consequence of (\ref{e.2.13}) is that
\begin{equation}
\label{e.2.14}
\| \varphi\|_{H^*\otimes\mathbf{C}}^2
	\le \|\varphi\|_{W^*\otimes\mathbf{C}}^2 \int_{W}\|w\| _{W}^{2}d\mu(w)
	= C_2\|\varphi\|_{W^*\otimes\mathbf{C}}^2.
\end{equation}

\subsection{Extensions of Lie algebras}
\label{s.ext}
Suppose $\mathfrak{v}$ is a Lie
algebra and $\mathrm{Der}(\mathfrak{v})$ is the set of derivations
on $\mathfrak{v}$.  That is, $\mathrm{Der}(\mathfrak{v})$ consists of all
linear maps $\rho:\mathfrak{v}\rightarrow\mathfrak{v}$ satisfying Leibniz's
rule:
\[ \rho([X,Y]_\mathfrak{v}) = [\rho(X),Y]_\mathfrak{v} +
	[X,\rho(Y)]_\mathfrak{v}. \]
$\mathrm{Der}(\mathfrak{v})$ forms a Lie algebra with Lie bracket defined by
the commutator:
\[ [\rho_1,\rho_2] = \rho_1\rho_2 - \rho_2\rho_1, \qquad \text{for }
	\rho_1,\rho_2\in\mathrm{Der}(\mathfrak{v}). \]
$\mathrm{Der}(\mathfrak{v})$ is a subset of linear maps on $\mathfrak{v}$,
so if $\mathfrak{v}$ is a normed vector space, one may
equip $\mathrm{Der}(\mathfrak{v})$ with the usual norm
\begin{equation}
\label{e.dernorm} 
\|\rho\|_0 = \sup\{\|\rho(X)\|_\mathfrak{v}:\|X\|_\mathfrak{v}=1\}.
\end{equation}
 
Now suppose that $\mathfrak{h}$ and $\mathfrak{v}$ are Lie algebras, and 
that there is a linear mapping 
\[ \alpha: \mathfrak{h} \rightarrow \mathrm{Der}(\mathfrak{v}) \]
and a skew-symmetric bilinear mapping
\[ \omega:\mathfrak{h}\times \mathfrak{h}\rightarrow\mathfrak{v}, \]
satisfying, for all $X,Y,Z\in \mathfrak{h}$,
\begin{equation}
\label{b1} 
\tag{B1}
[\alpha_X,\alpha_Y] - \alpha_{[X,Y]_\mathfrak{h}} 
	= \mathrm{ad}_{\omega(X,Y)} 
\end{equation}
and
\begin{equation}
\label{b2} 
\tag{B2}
\sum_{\text{cyclic}}\left(\alpha_X \omega(Y,Z) -
	\omega([X,Y]_\mathfrak{h},Z)\right)= 0.  
\end{equation}
Then, one may verify that, for $X_1+V_1,X_2+V_2\in \mathfrak{h}\oplus
\mathfrak{v} $,
\[ [X_1+V_1, X_2+V_2]_\mathfrak{g} :=  [X_1,X_2]_\mathfrak{h} + \omega(X_1,X_2) 
	+ \alpha_{X_1}V_2 -\alpha_{X_2}V_1 + [V_1,V_2]_\mathfrak{v} \]
defines a Lie bracket on $\mathfrak{g}:=\mathfrak{h}\oplus\mathfrak{v}$, 
and we say $\mathfrak{g}$ is an extension of $\mathfrak{h}$ over 
$\mathfrak{v}$. 
That is, $\mathfrak{g}$ is the Lie algebra with ideal
$\mathfrak{v}$ and quotient algebra
$\mathfrak{g}/\mathfrak{v}=\mathfrak{h}$.  The associated exact sequence is
\[ 0\rightarrow\mathfrak{v}\overset{\iota_1}{\longrightarrow}\mathfrak{g}
	\overset{\pi_2}{\longrightarrow}\mathfrak{h}\rightarrow0,\] 
where $\iota_1$ is
inclusion and $\pi_2$ is projection.  In fact, the
following theorem (see, for example, \cite{AMR00}) states that these are the 
only extensions of
$\mathfrak{h}$ over $\mathfrak{v}$.
\begin{thm}
Isomorphism classes of extensions of $\mathfrak{h}$ over $\mathfrak{v}$ (that
is, short exact sequences of Lie algebras $0
\rightarrow\mathfrak{v}\rightarrow\mathfrak{g}\rightarrow\mathfrak{h}
\rightarrow0$)
modulo the equivalence described by the commutative diagram of Lie algebra
homomorphisms
\begin{equation*}
\begin{CD}
0 @>>> \mathfrak{v}@>>> \mathfrak{g}@>>> \mathfrak{h}@>>> 0 \\
@.   @V{\mathrm{id}}VV   @V{\varphi}VV    @V{\mathrm{id}}VV \\
0 @>>> \mathfrak{v}@>>> \mathfrak{g}'@>>> \mathfrak{h}@>>> 0,
\end{CD}
\end{equation*}
correspond bijectively to equivalence classes of pairs of linear maps
$\alpha:\mathfrak{h} \rightarrow \mathrm{Der}(\mathfrak{v})$ and 
skew-symmetric bilinear maps $\omega:\mathfrak{h}\times
\mathfrak{h}\rightarrow\mathfrak{v}$
satisfying (\ref{b1}) and (\ref{b2}), where
$(\alpha,\omega)\equiv(\alpha',\omega')$ if there exists a linear
$b:\mathfrak{h}\rightarrow\mathfrak{v}$ such that
\[ \alpha_X' = \alpha_X + \mathrm{ad}_{b(X)}, \]
and
\begin{align*} 
\omega'(X,Y) &= \omega(X,Y) + \alpha_Xb(Y) - \alpha_Yb(X) - b([X,Y]) +
	[b(X),b(Y)]_\mathfrak{v}.
\end{align*}
The corresponding isomorphism $\varphi:\mathfrak{g}\rightarrow\mathfrak{g}'$ is
given by $\varphi(X+V)=X-b(X)+V$.
\end{thm}

When $\mathfrak{v}=V$ is an abelian Lie algebra, these pairs consist
of a Lie algebra homomorphism $\alpha:\mathfrak{h}\rightarrow gl(V)$ and
$\omega\in H^2(\mathfrak{h},V)$ is a Chevalley cohomology class with
coefficients in the $\mathfrak{h}$-module $V$ (see \cite{Fuks},
Chapter 1, Sections 3.1 and 4.5).  For definitions and details on
extensions of Lie algebras, see Section XIV.5 of \cite{CarEil99}.  Reference
\cite{AMR00} also gives a nice (although unpublished) summary.  Reference
\cite{Yankosky01} gives some conditions under which the extension of
$\mathfrak{h}$ over $\mathfrak{v}$ is nilpotent (when $\mathfrak{h}$
and $\mathfrak{v}$ are nilpotent); \cite{Sadetov05} gives a
characterization of extensions of a Lie algebra over a Heisenberg
Lie algebra.

\section{Semi-infinite Lie algebras and groups}
\label{s.semi-inf}
Throughout the rest of this paper $(W,H,\mu)$ will denote a real abstract
Wiener space,
and $\mathfrak{v}$ will denote a Lie algebra with $\mathrm{dim}(\mathfrak{v}) =
N<\infty$, equipped with an
inner product $\langle\cdot,\cdot\rangle_\mathfrak{v}$ and a continuous Lie
bracket $[\cdot,\cdot]_\mathfrak{v}$.  Note that this implies that there
exists a constant $c_0<\infty$ such that
\[ \|[X,Y]\|_\mathfrak{v} \le c_0 \|X\|_\mathfrak{v}\|Y\|_\mathfrak{v}, \]
for all $X,Y\in \mathfrak{v}$.  For simplicity,  
we will assume that $c_0\equiv1$.  Also, $\mathrm{Der}(\mathfrak{v})$ will 
denote the derivations of $\mathfrak{v}$, equipped with the norm defined 
in (\ref{e.dernorm}). 

\begin{defn}
Let $(W,H,\mu)$ be an abstract Wiener space and $\mathfrak{v}$ a finite
dimensional Lie algebra.  Then $\mathfrak{g}=W\oplus \mathfrak{v}$ endowed
with a Lie bracket satisfying
\begin{enumerate}
\item $[\mathfrak{g},\mathfrak{g} ] \subset \mathfrak{v}$, and
\item $[\cdot,\cdot]:\mathfrak{g}\times \mathfrak{g}
	\rightarrow \mathfrak{g}$ is continuous,
\end{enumerate}
will be called a {\em semi-infinite Lie algebra}.
\end{defn}

Motivated by the discussion in Section \ref{s.ext}, we may consider 
$W$ as an abelian Lie algebra and construct extensions of $W$ over
$\mathfrak{v}$.  So suppose there is a skew-symmetric
continuous bilinear mapping
\[ \omega:W\times W\rightarrow\mathfrak{v} \]
and a continuous linear mapping
\[ \alpha:W\rightarrow \mathrm{Der}(\mathfrak{v}) \]
such that $\alpha$ and $\omega$ satisfy (\ref{b1}) and (\ref{b2}), which in
this setting become
\begin{equation}
\label{c1} 
\tag{C1}[\alpha_X,\alpha_Y] = \mathrm{ad}_{\omega(X,Y)}
\end{equation}
and
\begin{equation}
\label{c2}
\tag{C2}
\alpha_X \omega(Y,Z) + \alpha_Y \omega(Z,X) + \alpha_Z \omega(X,Y) = 0, 
\end{equation}
for all $X,Y,Z\in W$.  Then we may define a Lie algebra structure
on $\mathfrak{g}:= W\oplus \mathfrak{v}$ via the Lie bracket
\[ [(X_1,V_1), (X_2,V_2)]_\mathfrak{g} := (0, \omega(X_1,X_2)
	+ \alpha_{X_1}V_2 - \alpha_{X_2}V_1 + [V_1,V_2]_\mathfrak{v}). \]
The vector space $\mathfrak{g}$ is also a Banach space in the norm
\[ \|(w,v)\|_{\mathfrak{g}} := \|w\|_W + \|v\|_\mathfrak{v}, \]
and $\mathfrak{g}_{CM}
:= H\oplus\mathfrak{v}$ is a Hilbert space with respect to the inner product
\[ \langle (A,a),(B,b)\rangle_{\mathfrak{g}_{CM}}
	:= \langle A,B\rangle_H + \langle a,b\rangle_\mathfrak{v}. \]
The associated Hilbertian norm on $\mathfrak{g}_{CM}$ is given by
\[ \|(A,a)\|_{\mathfrak{g}_{CM}} := \sqrt{\|A\|_H^2+\|a\|_\mathfrak{v}^2}. \]

\begin{notation}
\label{n.unorm}
Let 
\[ \|\omega\|_0 := \sup \{\|\omega(w_1,w_2)\|_\mathfrak{v} : \|w_1\|_W =
	\|w_2\|_W =1 \} \]
and
\[ \|\alpha\|_0 := \sup \{\|\alpha_w v\|_\mathfrak{v} : \|w\|_W =
	\|v\|_\mathfrak{v} = 1 \} \]
be the uniform norms of $\omega$ and $\alpha$, which are finite by their
assumed continuity.
\end{notation}
It will be useful to note that
\begin{equation}
\label{e.3} 
\|[\cdot,\cdot]\|_0 
	:= \sup\{ \|[g_1,g_2]\|_\mathfrak{v}: 
	\|g_1\|_\mathfrak{g} = \|g_2\|_\mathfrak{g} = 1 \}
	\le \|\omega\|_0 + 2\|\alpha\|_0 + 1 < \infty,
\end{equation}
and similarly
\begin{equation}
\label{e.Comega} 
C := C(\omega,\alpha)
	:= \sup \{\|[h,k]\|_\mathfrak{v} :
	\|h\|_{\mathfrak{g}_{CM}}=\|k\|_{\mathfrak{g}_{CM}}=1\}
	\le \|[\cdot,\cdot]\|_0 < \infty.
\end{equation}
Thus, for all $\ell=1,\cdots,r-1$,
\[ \|\mathrm{ad}^\ell_h k\|_\mathfrak{v} \le
	C^\ell\|h\|_{\mathfrak{g}_{CM}}^\ell \|k\|_{\mathfrak{g}_{CM}}. \]

If $\mathfrak{v}$ is nilpotent, $\omega$ and $\alpha$ may be chosen so 
that $\mathfrak{g}$ is a nilpotent Lie
algebra (see Section \ref{s.examples} for some examples).  
For $\mathfrak{g}$ nilpotent of step $r$, 
the Baker-Campbell-Hausdorff-Dynkin formula implies that
\[ \log(e^A e^B) = A+B+\sum_{k=1}^{r-1} 
		\sum_{(n,m)\in\mathcal{I}_k}  
		a_{n,m}^k\mathrm{ad}_A^{n_1} \mathrm{ad}_B^{m_1} \cdots
		\mathrm{ad}_A^{n_k} \mathrm{ad}_B^{m_k} A,
\]
for all $A,B\in\mathfrak{g}$, where 
\begin{equation} 
\label{e.ak}
a_{n,m}^k := \frac{(-1)^k}{(k+1)m!n!(|n|+1)}, 
\end{equation}
$\mathcal{I}_k := \{(n,m)\in\mathbb{Z}_+^k\times\mathbb{Z}_+^k : 
n_i+m_i>0 \text{ for all } 1\le i\le k \}$, and for each multi-index
$n\in\mathbb{Z}_+^k$,
\[ n!= n_1!\cdots n_k! \quad \text{ and } \quad |n|=n_1+\cdots+n_k, \]
see, for example, \cite{duiskolk}.  Since $\mathfrak{g}$ is nilpotent of step
$r$, 
\[ \mathrm{ad}_A^{n_1} \mathrm{ad}_B^{m_1} \cdots
		\mathrm{ad}_A^{n_k} \mathrm{ad}_B^{m_k} A = 0 \quad
\text{if } |n|+|m|\ge r. \]
for $A,B\in\mathfrak{g}$.  In particular, one may verify that
\begin{align}
\label{e.mult}
g\cdot h 
	&= g+h+\sum_{k=1}^{r-1} 
		\sum_{(n,m)\in\mathcal{I}_k}  
		a_{n,m}^k\mathrm{ad}_g^{n_1} \mathrm{ad}_h^{m_1} \cdots
		\mathrm{ad}_g^{n_k} \mathrm{ad}_h^{m_k} g
\end{align}
defines a group structure on $\mathfrak{g}$.  Note that $g^{-1}=-g$ and 
the identity $\mathbf{e}=(0,0)$.

\begin{defn}
When we wish to emphasize the group structure on $\mathfrak{g}$, we will
denote $\mathfrak{g}$ by $G$.  Similarly, when we wish to view
$\mathfrak{g}_{CM}$ as a subgroup of $G$, it will be denoted by $G_{CM}$ and
will be called the {\em Cameron-Martin subgroup}.
\end{defn}
(Since $\mathfrak{g}$ is simply connected and nilpotent, 
the exponential map is a global diffeomorphism (see, for
example, Theorems 3.6.2 of \cite{Varadarajan} or 1.2.1 of \cite{CorGrn90}), and we may identify
$\mathfrak{g}$ and $G$ under exponential coordinates.  In particular,
we may view $\mathfrak{g}$ as both a Lie algebra and Lie group.)

\begin{lem}
\label{l.cts}
The Banach space topologies on $\mathfrak{g}$ 
and $\mathfrak{g}_{CM}$ make $G$ and $G_{CM}$ into topological groups.
\end{lem}
\begin{proof}
Since $\mathfrak{g}$ and $\mathfrak{g}_{CM}$ are topological vector spaces, 
$g\mapsto g^{-1}=-g$ and $(g_1,g_2)\mapsto g_1+g_2$ are continuous by
definition.  The map $(g_1,g_2)\mapsto [g_1,g_2]$ is continuous in both 
the $\mathfrak{g}$ and $\mathfrak{g}_{CM}$ topologies by the estimates in
equations (\ref{e.3}) and (\ref{e.Comega}).
It then follows from (\ref{e.mult}) that $(g_1,g_2)\mapsto g_1\cdot g_2$
is continuous as well.
\end{proof}

\subsection{Examples}
\label{s.examples}

In this section, we give a few simple examples of semi-infinite Lie algebras.

\begin{example}
\label{ex-heis}
If $\mathfrak{v}$ is a finite dimensional inner product space, 
we may consider $\mathfrak{v}$ as an abelian Lie algebra, and
taking $\alpha\equiv0$ yields the infinite dimensional (step 2, stratified)
Heisenberg like Lie algebras described in \cite{DG08-2}.  
\end{example}

\begin{example}
\label{ex-beta}
Suppose $\mathfrak{v}$ is an $N$-dimensional nilpotent Lie algebra.  One
standard way to construct Lie algebra extensions is as follows.
Let $\beta:W\rightarrow \mathfrak{v}$ be a continuous linear map, and define
$\alpha:W\rightarrow\mathrm{Der}(\mathfrak{v})$ as the inner derivation
$\alpha_X:=\mathrm{ad}_{\beta(X)}$.  In this case, (\ref{c1}) and (\ref{c2})
are both satisfied if $\omega:W\times W\rightarrow \mathfrak{v}$
is given by $\omega(X,Y) := [\beta(X),\beta(Y)]_\mathfrak{v}$. Thus,
$\mathfrak{g}$ has Lie bracket
\[ [(X,V),(Y,U)]_\mathfrak{g} = (0,[\beta(X),\beta(Y)]_\mathfrak{v} 
	+ [\beta(X),U]_\mathfrak{v} - [\beta(Y),V]_\mathfrak{v} 
	+ [V,U]_\mathfrak{v}), \]
and, if $\mathfrak{v}$ is nilpotent Lie algebra of step $r$, then
$\mathfrak{g}$ is nilpotent of step $r$. 

One should note for this construction that, 
since $\beta$ is linear, we have the decomposition 
$W = \mathrm{Nul}(\beta) \oplus \mathrm{Nul}(\beta)^\perp$, where
$\mathrm{dim}(\mathrm{Nul}(\beta)^\perp) \le \mathrm{dim}(\mathfrak{v})=N$.
Thus, for $X=X_1+X_2, Y=Y_1+Y_2\in W$,
\[ \omega( X_1+X_2,Y_1+Y_2) = [ \beta(X_1+X_2),\beta(Y_1+Y_2)]
	= [\beta(X_2),\beta(Y_2)], \]
and $\omega$ is a map on $\mathrm{Nul}(\beta)^\perp\times
\mathrm{Nul}(\beta)^\perp$.  
Thus, $[\mathrm{Nul}(\beta),\mathrm{Nul}(\beta)]=\{0\}$ and similarly
$[\mathrm{Nul}(\beta),\mathfrak{v}]=\{0\}$.  So
\[ \mathfrak{g}= W\oplus\mathfrak{v} = \mathrm{Nul}(\beta) \oplus
	\mathrm{Nul}(\beta)^\perp\oplus\mathfrak{v} \]
is in a sense just an extension of the finite dimensional subspace
$\mathrm{Nul}(\beta)^\perp$ by $\mathfrak{v}$.
\end{example}

\begin{example}
\label{ex-beta2}
One can generalize the previous example by taking a linear map
$\beta:W\rightarrow\mathfrak{h}$, where $\mathfrak{h}$ is nilpotent Lie
algebra, and constructing an extension of $\mathfrak{h}$ by a nilpotent Lie
algebra.  For the sake of a concrete example, consider the following.

Let 
\[ W=W(\mathbb{R}^3)=\{\sigma:[0,1]\rightarrow\mathbb{R}^3 : \sigma \text{ is
continuous and } \sigma(0)=0\} \]
and 
\[ H=\left\{\sigma\in W : \sigma \text{ is
absolutely continuous and } \int_0^1 \|\dot{\sigma}(s)\|^2\,ds < \infty
	\right\}, \]
so that $(W,H)$ is standard Wiener space.  Let
$\mathfrak{v}=\mathbb{R}^3$ be an abelian Lie algebra. 
Let $\bar{\sigma}= \int_0^1 \sigma(s)\,ds = 
(\bar{\sigma}_1,\bar{\sigma}_2,\bar{\sigma}_3)$, and
define $\omega:W\times W\rightarrow \mathbb{R}^3$ by
\[ \omega(\sigma,\tau) 
	= \left(\bar{\sigma}_1\bar{\tau}_2 - \bar{\tau}_1\bar{\sigma}_2,
		\bar{\sigma}_2\bar{\tau}_3 - \bar{\tau}_2\bar{\sigma}_3,
		0 \right) \]
and $\alpha_\sigma:\mathbb{R}^3\rightarrow \mathbb{R}^3$ by
\[ \alpha_\sigma (x,y,z) = (0,0,\bar{\sigma}_1y-\bar{\sigma}_3x). \]
Then $\alpha_\sigma\alpha_\tau=0$ and (\ref{c1}) is trivially satisfied.
Using that
\[ \alpha_\kappa \omega(\sigma,\tau) = (0,0, \bar{\kappa}_1
	(\bar{\sigma}_2\bar{\tau}_3 - \bar{\tau}_2\bar{\sigma}_3)
	-\bar{\kappa}_3(\bar{\sigma}_1\bar{\tau}_2 -
	 \bar{\tau}_1\bar{\sigma}_2 )) \]
one may verify that (\ref{c2}) is satisfied.  Thus, the Lie bracket for this
extension $\mathfrak{g}=W\oplus\mathbb{R}^3$ is given by
\[ [(\sigma,v),(\tau,u)]
	= (0,\bar{\sigma}_1\bar{\tau}_2 - \bar{\tau}_1\bar{\sigma}_2,
		\bar{\sigma}_2\bar{\tau}_3 - \bar{\tau}_2\bar{\sigma}_3, 
		\bar{\sigma}_1u_2-\bar{\sigma}_3u_1 +
		\bar{\tau}_1v_2-\bar{\tau}_3v_1), \]
\[ [(\kappa,w),[(\sigma,v),(\tau,u)]]
	= (0,0,0, \bar{\kappa}_1
	(\bar{\sigma}_2\bar{\tau}_3 - \bar{\tau}_2\bar{\sigma}_3)
	-\bar{\kappa}_3(\bar{\sigma}_1\bar{\tau}_2 -
	 \bar{\tau}_1\bar{\sigma}_2)), \]
and all higher order brackets are 0.

Note that this construction corresponds to the extension
$\mathfrak{g}= \mathbb{R}^3\oplus\mathbb{R}^3$, the $4\times4$ upper 
triangular matrices.  To see this, let $U=\mathbb{R}^3$ and
$V=\mathfrak{v}=\mathbb{R}^3$, and define $\omega':U\times U\rightarrow V$ by
\begin{align*}
\omega( (a,b,c),&(a',b',c')) \\
	&=\begin{pmatrix} 0 & a & 0 & 0 \\
		0 & 0 & b & 0\\
		0 & 0 & 0 & c \\
		0 & 0 & 0 & 0 \end{pmatrix} 
		\begin{pmatrix} 0 & a' & 0 & 0 \\
		0 & 0 & b' & 0\\
		0 & 0 & 0 & c' \\
		0 & 0 & 0 & 0 \end{pmatrix}  
		-\begin{pmatrix} 0 & a' & 0 & 0 \\
		0 & 0 & b' & 0\\
		0 & 0 & 0 & c' \\
		0 & 0 & 0 & 0 \end{pmatrix} 
		\begin{pmatrix} 0 & a & 0 & 0 \\
		0 & 0 & b & 0\\
		0 & 0 & 0 & c \\
		0 & 0 & 0 & 0 \end{pmatrix} \\
	&= \begin{pmatrix} 0 & 0 & ab'- ba' & 0 \\
		0 & 0 & 0 & bc'-cb' \\
		0 & 0 & 0 & 0 \\
		0 & 0 & 0 & 0 \end{pmatrix},
\end{align*}
and $\alpha':U\rightarrow gl(V)$ by
\begin{align*} 
\alpha_{(a,b,c)}&(x,y,z) \\
	&= \begin{pmatrix} 0 & a & 0 & 0 \\
		0 & 0 & b & 0\\
		0 & 0 & 0 & c \\
		0 & 0 & 0 & 0 \end{pmatrix} 
	\begin{pmatrix} 0 & 0 & x & z \\
		0 & 0 & 0 & y\\
		0 & 0 & 0 & 0 \\
		0 & 0 & 0 & 0 \end{pmatrix} 
	- \begin{pmatrix} 0 & 0 & x & z \\
		0 & 0 & 0 & y\\
		0 & 0 & 0 & 0 \\
		0 & 0 & 0 & 0 \end{pmatrix} 
	\begin{pmatrix} 0 & a & 0 & 0 \\
		0 & 0 & b & 0\\
		0 & 0 & 0 & c \\
		0 & 0 & 0 & 0 \end{pmatrix} \\
	&= \begin{pmatrix} 0 & 0 & 0 & ay-cx \\
		0 & 0 & 0 & 0\\
		0 & 0 & 0 & 0 \\
		0 & 0 & 0 & 0 \end{pmatrix}.
\end{align*}
Then $\omega=\omega'\circ\beta$ and $\alpha=\alpha'\circ\beta$ where
$\beta:W\rightarrow U$ is given by
$\beta(\sigma)=(\bar{\sigma}_1,\bar{\sigma}_2,\bar{\sigma}_3)$.
\end{example}

\begin{example}
Consider
$\mathfrak{v}=\mathbb{R}^n\oplus\mathbb{R}$ as an abelian Lie algebra.  
For $\omega:W\times W\rightarrow\mathbb{R}^n$, 
we may write $\omega = (\omega_1,\cdots,\omega_n)$, where $\omega_i:W\times
W\rightarrow\mathbb{R}$ are bilinear, anti-symmetric, continuous maps.
Similarly, for $\alpha:W\times\mathbb{R}^n\rightarrow\mathbb{R}$, we have
$\alpha_i(\cdot) = \alpha_\cdot e_i$, where $\{e_i\}_{i=1}^n$ is the 
standard basis for $\mathbb{R}^n$.  Thus,
\[ \alpha_w (a_1,\ldots,a_n) = \sum_{i=1}^n a_i \alpha_i(w). \]
Then $\alpha$ and $\omega$ satisfy (\ref{c2}) as long as
\[ \alpha_1\wedge\omega_1+\cdots+\alpha_n\wedge\omega_n=0. \]

In the case $n=1$, this is not very interesting, since $
\alpha\wedge\omega=0$ implies that $\omega=\alpha\wedge\beta$ for some 
$\beta\in W^*$.

For $n=2$, we have $\mathfrak{v}=\mathbb{R}^2\oplus 
\mathbb{R}$.  Let 
$\Omega:W\times W\rightarrow\mathbb{R}$ be bilinear, antisymmetric, and
continuous,
and $\gamma:W\rightarrow\mathbb{R}$ be linear and continuous.  Then define
$\omega:W\times W\rightarrow\mathbb{R}^2$ by $\omega=(\Omega,\Omega)$
and $\alpha:W\times\mathbb{R}^2\rightarrow\mathbb{R}$ by $\alpha_1=
\gamma$ and $\alpha_2=-\gamma$, so that, for any $u,w\in W$ and
$v=(v_1,v_2)\in\mathbb{R}^2$,
\[ \omega(w,u) = (\Omega(w,u),\Omega(w,u)) \quad \text{ and } 
	\quad \alpha_w v = \gamma(w)(v_1-v_2). \]
Note that, for any $w,u,h\in W$, $\omega$ and $\alpha$ satisfy
\[ \alpha_h\omega(w,u) = \alpha_h(\Omega(w,u),\Omega(w,u))
	= \gamma(h)(\Omega(w,u)-\Omega(w,u))= 0. \]
Thus, for any $(w,v,x),(w',v',x'),(w'',v'',x'')\in W\oplus\mathfrak{v}$,
\begin{align*} 
[(w,v,x),(w',v',x')] 
	&= (0,\omega(w,w'),\alpha_wv'-\alpha_{w'}v) \\
	&= (0,(\Omega(w,w'),\Omega(w,w')),\gamma(w)(v_1'-v_2')
			+ \gamma(w')(v_1-v_2)), 
\end{align*}
\[ [(w'',v'',x''),[(w,v,x),(w',v',x')]]
	= (0,0,\alpha_{w''}\omega(w,w')) = 0, \]
and $\mathfrak{g}$ is a step 2 Lie algebra.  The group operation is given by
\begin{align*} 
(w,v,x)\cdot(w',v',x')
	&= (w+w',v+v'+\frac{1}{2}(\Omega(w,w'),\Omega(w,w')), \\
	&\qquad\quad	x+x'+\frac{1}{2}(\gamma(w)(v_1'-v_2')
			+ \gamma(w')(v_1-v_2)).
\end{align*}

As an example of a particular appropriate $\Omega$ and $\gamma$, again let
$W=W(\mathbb{R}^3)$ and $H$ be as in Example \ref{ex-beta2}.  Suppose 
$\varphi$ is an anti-symmetric bilinear form on $\mathbb{R}^3$,
$\rho:\mathbb{R}^3\rightarrow\mathbb{R}$ is a linear map, and let $\eta$
be finite measure on $[0,1]$.  Then we may define
\[ \Omega(\sigma, \tau) = \int_0^1 \varphi(\sigma(s),\tau(s))\,d\eta(s) \]
and
\[ \gamma(\sigma) = \int_0^1 \rho(\sigma(s)) \,d\eta(s). \]
\end{example}


\begin{example}
Here we make a slight modification on the previous example to construct a
stratified step 3 Lie algebra.  Let $\mathfrak{v}=\mathbb{R}^6
= \mathbb{R}^3\oplus\mathbb{R}^2\oplus\mathbb{R}$ be an abelian Lie algebra.
Let $\Omega$ and $\gamma$ be as in the previous example.
Define $\omega:W\times W\rightarrow\mathbb{R}^3$ by
\[ \omega(w,u) = (\Omega(w,u),\Omega(w,u),\Omega(w,u)) \]
and $\alpha:W\times\mathfrak{v}\rightarrow\mathfrak{v}$ by
\[ \alpha_w ( (v_1,v_2,v_3), (x_1,x_2),y)
	= (0, (\gamma(w)(v_1-v_2),\gamma(w)(v_2-v_3)),\gamma(w)(x_1-x_2)) \]
(so $\alpha_w$ is a particular element of the $6\times6$ strictly 
lower triangular matrices).
Then $\alpha_w\alpha_u=\alpha_u\alpha_w$ and so $\alpha$ satisfies (\ref{c1}),
and also
\[ \alpha_v\omega(w,u) = (0,(\gamma(v)(\Omega(w,u)-\Omega(w,u)),
	\gamma(v)(\Omega(w,u)-\Omega(w,u))),0) = 0,  \]
so $\alpha$ and $\omega$ satisfy (\ref{c2}) trivially.  The Lie bracket is
given by
\[ 
[(w,v,x,y),(w',v',x',y')]
	= (0,\omega(w,w'),\alpha_w v'- \alpha_{w'} v, 
		\alpha_w x' - \alpha_{w'}x), \]
or, more explicitly, this may be written componentwise as
\[ [(w,v,x,y),(w',v',x',y')]_2
	= (\Omega(w,w'),\Omega(w,w'),\Omega(w,w'))\in\mathbb{R}^3, \]
\begin{multline*}
[(w,v,x,y),(w',v',x',y')]_3 \\
	= (\gamma(w)(v_1'-v_2')-
		\gamma(w')(v_1-v_2), \gamma(w)(v_2'-v_3')-\gamma(w')(v_2-v_3))
		\in\mathbb{R}^2, 
\end{multline*}
and
\[ [(w,v,x,y),(w',v',x',y')]_4
	= \gamma(w)(x_1'-x_2') - \gamma(w')(x_1-x_2)\in\mathbb{R}.
\]
Thus,
\begin{align*} 
[(w'',v'',x'',y''),&[(w,v,x,y),(w',v',x',y')]] \\
	&= (0,0,\alpha_{w''}\omega(w,w'), 
		\alpha_{w''}(\alpha_w v' - \alpha_{w'} v)) \\
	&= (0,0,0, \alpha_{w''}\alpha_w v' -\alpha_{w''} \alpha_{w'} v) \\
	&= (0,0,0, \gamma(w'')\gamma(w)(v_1'-v_3') -
			\gamma(w'')\gamma(w')(v_1-v_3)),
\end{align*}
and all higher order brackets are 0.  So for $g=(w,v,x,y)$ and $g'=
(w',v',x',y')$, the group operation is given by
\begin{align*}
(g\cdot g')_1 &= w+w' \\
(g\cdot g')_2 &= v+v'+ \frac{1}{2}\omega(w,w') \\
(g\cdot g')_3 &= x+x' + \frac{1}{2} (\alpha_w v' -\alpha_{w'}v) \\
(g\cdot g')_4 &= y+y'+ \frac{1}{2} (\alpha_wx'-\alpha_{w'}x) 
		+ \frac{1}{12} (\alpha_w^2v'+ \alpha_{w'}^2 v 
		- \alpha_w\alpha_{w'} (v-v')).
\end{align*}
Clearly, this example may be further modified to make
nilpotent Lie algebras of arbitrary step.

\end{example}

\subsection{Hilbert-Schmidt norms}
\label{s.norm}

In this section, we will show that the assumed continuity of $\omega$ and
$\alpha$ makes the Lie bracket into a Hilbert-Schmidt operator on
$\mathfrak{g}_{CM}$.  This result will be needed later in guaranteeing that
our stochastic integrals are well-defined.

\begin{notation}
Let $H_1,\ldots,H_n$ and $V$ be Hilbert spaces, and let 
$\{h_j^i\}_{j=1}^{\mathrm{dim}(H_i)}$ denote an orthonormal basis 
for each $H_i$.  If $\rho:H_1\times\cdots\times H_n\rightarrow V$ is a 
multilinear map, then the Hilbert-Schmidt norm of $\rho$ is defined by
\[ \|\rho\|_2^2 : = \|\rho\|_{H_1^*\otimes\cdots\otimes H_n^*\otimes V}
	= \sum_{j_1,\ldots,j_n} \|\rho(h_{j_1}^1,\ldots,h_{j_n}^n)\|^2_V. \]
In particular, for $H$ an infinite dimensional Hilbert space with orthonormal
basis $\{h_i\}_{i=1}^\infty$,  
$\rho:H^{\otimes n}\rightarrow V$ is Hilbert-Schmidt if
\[ \|\rho\|_2^2 = \|\rho\|_{(H^*)^{\otimes n}\otimes V}
	= \sum_{j_1,\ldots,j_n=1}^\infty \|\rho(h_{j_1},\ldots,h_{j_n})\|^2_V
	<\infty. \]
One may verify directly that these norms are independent of 
the chosen bases.
\end{notation}

\begin{prop}
\label{p.normest}
For all $w\in W$ and $x\in\mathfrak{v}$,
\begin{equation}
\label{e.alphaest}
\|\alpha_w \cdot\|^2_{\mathfrak{v}^*\otimes\mathfrak{v}} 
	\le N\|\alpha\|^2_0 \|w\|^2_W
\quad \text{ and } \quad
\|\alpha_\cdot x\|^2_{H^*\otimes\mathfrak{v}}
	\le C_2\|\alpha\|_0^2\|x\|_\mathfrak{v}^2,
\end{equation}
where $C_2$ is as in equation (\ref{e.2.2}).
Also,
\begin{equation}
\label{e.omegaest}
\|\omega(w,\cdot)\|^2_{H^*\otimes\mathfrak{v}} 
	\le C_2\|\omega\|_0^2 \|w\|_W^2. 
\end{equation}
Furthermore,
\[ \|\alpha\|^2_2
	\le NC_2 \|\alpha\|_0^2 < \infty \quad \text{ and }
	\quad \|\omega\|_2^2 
	\le C_2^2 \|\omega\|_0^2 < \infty. \]
\end{prop}
\begin{proof}
Let $\{e_i\}_{i=1}^N$ be an orthonormal basis of $\mathfrak{v}$.  Then, for
any $w\in W$,
\[
\|\alpha_w \cdot\|^2_{\mathfrak{v}^*\otimes\mathfrak{v}}
	= \sum_{i=1}^N \|\alpha_w e_i\|_\mathfrak{v}^2
	\le \sum_{i=1}^N \|\alpha\|^2_0\|w\|^2_W \|e_i\|_\mathfrak{v}^2
	= N \|\alpha\|_0^2 \|w\|^2_W.
\]
For fixed $x\in\mathfrak{v}$, $\alpha_\cdot x:W\rightarrow\mathfrak{v}$ is a
continuous linear map.  Thus, equation (\ref{e.2.13}) gives
\begin{align*} 
\|\alpha_\cdot x\|^2_{H^*\otimes\mathfrak{v}}
	&= \int_W \|\alpha_w x\|^2_\mathfrak{v}\,d\mu(w) \\
	&\le \int_W \|\alpha\|_0^2\|w\|_W^2\|x\|_\mathfrak{v}^2 \,d\mu(w)
	= C_2\|\alpha\|_0^2\|x\|_\mathfrak{v}^2.
\end{align*}
Similarly, for fixed $w\in W$ and $\omega(w,\cdot):W\rightarrow\mathfrak{v}$,
\begin{align*}
\|\omega(w,\cdot)\|^2_{H^*\otimes\mathfrak{v}}
	&= \int_W \|\omega(w,w')\|^2_\mathfrak{v} \,d\mu(w') \\
	&\le \int_W  \|\omega\|_0^2 \|w\|_W^2 \|w'\|^2_W \,d\mu(w') 
	=  C_2\|\omega\|_0^2 \|w\|_W^2.
\end{align*}

Since $w\mapsto\alpha_w$ is a continuous linear map from $W$ to 
$\mathfrak{v}^*\otimes\mathfrak{v}$, it follows from equations (\ref{e.2.13}) 
and (\ref{e.2.14}) that
\begin{align*}
\|\alpha\|_2^2 
	= \int_W \|\alpha_w \cdot\|_{\mathfrak{v}^*\otimes\mathfrak{v}}^2
		\,d\mu(w) 
	&\le \int_W N \|\alpha\|_0^2\|w\|_W^2\,d\mu(w) = NC_2\|\alpha\|_0^2,
\end{align*}
and since $w\mapsto\omega(w,\cdot)$ is a continuous linear map from
$W$ to $H^*\otimes\mathfrak{v}$, 
\begin{align*}
\|\omega\|_2^2 
	&= \|h\mapsto\omega(h,\cdot)\|^2_{H^*\otimes(H^*\otimes\mathfrak{v})}
	= \int_W \|\omega(w,\cdot)\|_{H^*\otimes\mathfrak{v}}^2\,d\mu(w) \\
	&\le \int_W C_2 \|\omega\|_0^2\|w\|_W^2\,d\mu(w) = C_2^2\|\omega\|_0^2.
\end{align*}
\end{proof}

This proposition easily gives the following result. 

\begin{cor}
\label{c.normest}
For all $m\ge2$, $[[[\cdot,\cdot],\ldots],\cdot]:\mathfrak{g}_{CM}^{\otimes m}
\rightarrow \mathfrak{v}$ is Hilbert-Schmidt.
\end{cor}
\begin{proof}
For $m=2$, this follows from the previous proposition and the
continuity of the Lie bracket on $\mathfrak{v}$, since taking
$\{h_i\}_{i=1}^\infty=\{k_i\}_{i=1}^\infty\cup\{e_j\}_{j=1}^N$, where 
$\{k_i\}_{i=1}^\infty$ and $\{e_j\}_{j=1}^N$ are orthonormal bases of $H$ and
$\mathfrak{v}$, respectively, gives
\begin{align*}
\|[\cdot,\cdot]\|_2^2
	&=	\|[\cdot,\cdot]\|^2
		_{\mathfrak{g}_{CM}^*\otimes\mathfrak{g}_{CM}^*\otimes
		\mathfrak{v}} 
	= \sum_{i_1,i_2=1}^\infty 
		\|[ [h_{i_1},h_{i_2}]\|_\mathfrak{v}^2 \\
	&= \sum_{i_1,i_2=1}^\infty
		\|\omega(k_{i_1},k_{i_2})\|_\mathfrak{v}^2 
		+ \sum_{i_1=1}^\infty\sum_{j_2=1}^N
			\|\alpha_{k_{i_1}}e_{j_2}\|_\mathfrak{v}^2 \\
	&\qquad + \sum_{i_2=1}^\infty\sum_{j_1=1}^N
			\|\alpha_{k_{i_2}}e_{j_1}\|_\mathfrak{v}^2 
		+ \sum_{j_1,j_2=1}^N \|[e_{j_1},e_{j_2}] \|_\mathfrak{v}^2 \\
	&= \|\omega\|_2^2 + 2\|\alpha\|_2^2 + N <\infty.
\end{align*}
Now assume the statement is true for all $m=2,\ldots,\ell$.  
Consider $m=\ell+1$.
Writing $[[h_{i_1},h_{i_2}],\cdots,h_{i_\ell}]\in 
\mathfrak{v}$ in terms of the orthonormal basis
$\{e_j\}_{j=1}^N$ and using 
multiple applications of the Cauchy-Schwarz inequality gives
\begin{align*}
\|[[[\cdot,&\cdot],\ldots],\cdot]\|_2^2
	= \|[[[\cdot,\cdot],\ldots],\cdot]\|_{(\mathfrak{g}_{CM}^*)^{\otimes \ell+1}
		\otimes\mathfrak{v}} \\
	&= \sum_{i_1,\ldots,i_{\ell+1}=1}^\infty 
		\|[[[h_{i_1},h_{i_2}],\cdots,h_{i_\ell}],h_{i_{\ell+1}}]\|
		_\mathfrak{v}^2 \\
	&= \sum_{i_1,\ldots,i_{\ell+1}=1}^\infty \left\|\sum_{j=1}^N 
		[e_j,h_{i_{\ell+1}}]
		\langle e_j, [[h_{i_1},h_{i_2}],\cdots,h_{i_\ell}]\rangle
		\right\|_\mathfrak{v}^2 \\
	&\le N \sum_{i_1,\ldots,i_{\ell+1}=1}^\infty \sum_{j=1}^N 
		\| [e_j,h_{i_{\ell+1}}]\|_\mathfrak{v}^2 
		|\langle e_j,[[h_{i_1},h_{i_2}],\cdots,h_{i_\ell}]\rangle|^2 \\
	&\le N \left(\sum_{i_{\ell+1}=1}^\infty \sum_{j=1}^N 
			\| [e_j,h_{i_{\ell+1}}] \|^2_\mathfrak{v} \right)
		\left(\sum_{i_1,\ldots,i_\ell=1}^\infty \sum_{j=1}^N 
		|\langle e_j,[[h_{i_1},h_{i_2}],\cdots,h_{i_\ell}]\rangle|^2\right) \\
	&\le N \|[\cdot,\cdot]\|^2_{\mathfrak{g}_{CM}^{\otimes 2}\otimes 
		\mathfrak{v}} \cdot
		\|[[[\cdot,\cdot],\ldots],\cdot]\|^2
		_{\mathfrak{g}_{CM}^{\otimes \ell}\otimes \mathfrak{v}},
\end{align*}
where in the penultimate inequality we have used that all terms in the sums
are positive.  The last line is finite by the induction hypothesis.
\end{proof}

\subsection{Length and distance}
\label{s.length}

In this section, we define the Riemannian distance on $G_{CM}$ 
and show that the topology induced by this metric is equivalent to the
Hilbert topology induced by $\|\cdot\|_{\mathfrak{g}_{CM}}$.

For $g\in G,$ let $L_g:G\rightarrow G$ and $R_g:G\rightarrow G$
denote left and right multiplication by $g$, respectively.  As $G$
is a vector space, to each $g\in G$ we can associate the tangent
space $T_g G$ to $G$ at $g$, which is naturally isomorphic to $G$.

\begin{notation}
For $f:G\rightarrow\mathbb{R}$ a Frech\'{e}t smooth function and  
$v,x\in G$ and $h\in\mathfrak{g}$, let
\[ f'(x)h := \partial_h f(x) = \frac{d}{dt}\bigg|_{0}f(x+th), \]
and let $v_x \in T_x G$ denote the tangent vector
satisfying $v_xf=f'(x)v$.  If $\sigma(t)$ is any smooth curve in
$G$ such that $\sigma(0) = x$ and $\dot{\sigma}(0)=v$ (for example,
$\sigma(t) = x+tv$), then
\[ L_{g*} v_x = \frac{d}{dt}\bigg|_0 g\cdot \sigma(t). \]
\end{notation}

\begin{notation}
\label{n.length}
Let $T>0$ and $C^1([0,T],G_{CM})$ denote the collection of $C^1$-paths
$g:[0,T]\rightarrow G_{CM}$.  The length of $g$ is defined as 
\[ \ell_{CM}(g) := \int_0^T \|L_{g^{-1}(s)*}g'(s)\|_{\mathfrak{g}_{CM}}\,ds.
\]
The Riemannian distance between $x,y\in G_{CM}$ then takes the usual form
\[ d_{CM}(x,y) := \inf\{\ell_{CM}(g): g\in C^1([0,T],G_{CM}) \text{ such that }
	g(0)=x \text{ and } g(T)=y \}. \]
Note that the value of $T$ in the definition of $d_{CM}$ is irrelevant
since the length functional is invariant under reparameterization.
\end{notation}

\begin{prop}
For $g,x\in G$ and $v_x\in T_x G$,
\begin{equation}
\label{e.gender}
\begin{split}
L_{g*} v_x &=  v + \sum_{k=1}^{r-1} 
		\sum_{(n,m)\in\mathcal{I}_k} a_{n,m}^k \times \\
	& \sum_{\tiny \begin{array}{cc}j\in\{1,\ldots,k\} \\
		m_j>0\end{array}} \sum_{\ell=0}^{m_j-1} 
		\mathrm{ad}_g^{n_1} \mathrm{ad}_x^{m_1} \cdots
		\mathrm{ad}_g^{n_j} \mathrm{ad}_x^{\ell}  \mathrm{ad}_v
		\mathrm{ad}_x^{m_j-\ell-1}\mathrm{ad}_g^{n_{j-1}} 
		\cdots  \mathrm{ad}_g^{n_k} \mathrm{ad}_x^{m_k} g,
\end{split}
\end{equation}
where $a_{n,m}^k$ are the coefficients in the group multiplication given
in equation (\ref{e.ak}).
\end{prop}
\begin{proof}
The proof is a simple computation.  Let $x(t)=x+tv$, and first note that
\begin{multline*}
\frac{d}{dt}\bigg|_0 \mathrm{ad}_g^{n_1} \mathrm{ad}_{x(t)}^{m_1} \cdots
		\mathrm{ad}_g^{n_k} \mathrm{ad}_{x(t)}^{m_k} g \\
	= \sum_{\tiny \begin{array}{cc}j\in\{1,\ldots,k\} \\
		m_j>0\end{array}} \sum_{\ell=0}^{m_j-1} 
		\mathrm{ad}_g^{n_1} \mathrm{ad}_x^{m_1} \cdots
		\mathrm{ad}_g^{n_j} \mathrm{ad}_x^{\ell}  \mathrm{ad}_v
		\mathrm{ad}_x^{m_j-\ell-1}\mathrm{ad}_g^{n_{j-1}}
		\cdots \mathrm{ad}_g^{n_k} \mathrm{ad}_x^{m_k} g.
\end{multline*}
Then using (\ref{e.mult}) and plugging this into
\begin{align*} 
&L_{g*} v_x 
	= \frac{d}{dt}\bigg|_0 g\cdot x(t) \\
	&\quad= \frac{d}{dt}\bigg|_0 \left( g + x(t)
		+ \sum_{k=1}^{r-1} 
		\sum_{(n,m)\in\mathcal{I}_k} a_{n,m}^k \mathrm{ad}_g^{n_1} \mathrm{ad}_{x(t)}^{m_1} \cdots
		\mathrm{ad}_g^{n_k} \mathrm{ad}_{x(t)}^{m_k} g \right) 
\end{align*}
yields the desired result.
\end{proof}

\begin{example}[The step 3 case]
When $r=3$, the group operation is
\[ g\cdot h = g+h +\frac{1}{2}[g,h] + \frac{1}{12}([g,[g,h]] + [h,[h,g]]). \]
Thus,
\begin{align*}
L_{g*}v_x &= \frac{d}{dt}\bigg|_0 g\cdot x(t) \\
	&= \frac{d}{dt}\bigg|_0 \left( g+x(t) +\frac{1}{2}[g,x(t)] 
		+ \frac{1}{12}([g,[g,x(t)]] + [x(t),[x(t),g]])\right) \\
	&= v + \frac{1}{2}[g,v] + \frac{1}{12}([g,[g,v]]+
		[v,[x,g]] + [x,[v,g]]).
\end{align*}
\end{example}

\begin{prop}
\label{p.length1}
There exists $K_1=K_1(a\wedge b)<\infty$ (for $a,b\ge0$) such that 
$K_1(0)=0$ and, for all $x,y\in G_{CM}$, 
\[
d_{CM}(x,y) \le 
	(1+K_1(\|x\|_{\mathfrak{g}_{CM}}\wedge\|y\|_{\mathfrak{g}_{CM}}))
	\|y-x\|_{\mathfrak{g}_{CM}} + o\left(\|y-x\|_{\mathfrak{g}_{CM}}^2\right).
\]
\end{prop}

\begin{proof}
For notational simplicity, let $T=1$.  If $g(s)$ is a
path in $C_{CM}^1$ for $0\le s\le1$, then, by equation (\ref{e.gender}),
taking $g=g^{-1}(s)$, $x=g(s)$, and $v_{g(s)}=g'(s)$,
\begin{align}
\notag
\ell_{CM}(g) 
	&= \int_0^1 \bigg\|g'(s) + \sum_{k=1}^{r-1} \sum_{(n,m)\in\mathcal{I}_k}
		a_{n,m}^k \sum_{m_j>0}\sum_{\ell=0}^{m_j-1}
		\mathrm{ad}_{g^{-1}(s)}^{n_1} \mathrm{ad}_{g(s)}^{m_1} 
		\\
	&\notag\qquad\qquad 
		\cdots \mathrm{ad}_{g^{-1}(s)}^{n_j}\mathrm{ad}_{g(s)}^{\ell} 
		\mathrm{ad}_{g'(s)} \mathrm{ad}_{g(s)}^{m_j-\ell-1} \cdots
		\mathrm{ad}_{g^{-1}(s)}^{n_k} \mathrm{ad}_{g(s)}^{m_k} g^{-1}(s)
		\bigg\|_{\mathfrak{g}_{CM}}\,ds \\
	&\notag
	= \int_0^1 \left\|g'(s) + \sum_{k=1}^{r-1} \sum_{(n,m)\in\mathcal{I}_k}
		(-1)^{|n|}1_{\{m_k>0\}} a_{n,m}^k		
		\mathrm{ad}_{g(s)}^{|m|+|n|} g'(s)
		\right\|_{\mathfrak{g}_{CM}}\,ds \\
	&\label{e.lg}
	= \int_0^1 \left\|g'(s) + \sum_{\ell=1}^{r-1}
		d_\ell \mathrm{ad}_{g(s)}^\ell g'(s)
		\right\|_{\mathfrak{g}_{CM}}\,ds, 
\end{align}
where 
\begin{equation} 
\label{e.dl}
d_\ell 
	:= \sum_{k=1}^\ell 
		\sum_{\tiny \begin{array}{cc}(n,m)\in\mathcal{I}_k\\
		|m|+|n|=\ell\end{array}}
		(-1)^{|n|}1_{\{m_k>0\}} a_{n,m}^k.
\end{equation}
Taking $g(s)=x+s(y-x)$ for $0\le s\le 1$, this gives
\begin{align*}
d&_{CM}(x,y) \le \ell_{CM}(g) \\
	&= \int_0^1 \left\|(y-x) + \sum_{\ell=1}^{r-1}
		d_\ell \mathrm{ad}_{x+s(y-x)}^\ell (y-x)
		\right\|_{\mathfrak{g}_{CM}} ds \\
	&= \int_0^1 \Bigg\|(y-x) + \sum_{\ell=1}^{r-1}
		d_\ell \sum_{\tiny \begin{array}{cc}
		(n,m)\in\mathcal{I}_\ell \\ |m|+|n|=\ell
		\end{array}} s^{|n|} 
		\mathrm{ad}_x^{m_1} \mathrm{ad}_{y-x}^{n_1}
		\cdots \mathrm{ad}_x^{m_\ell} \mathrm{ad}_{y-x}^{n_\ell} (y-x)
		\Bigg\|_{\mathfrak{g}_{CM}} ds.
\end{align*}
Splitting off all terms in the sum of order two or higher and evaluating the
integral gives
\begin{align*}
d_{CM}(x,y)
	&\le \Bigg\|(y-x) 
		+ \sum_{\ell=1}^{r-1} d_\ell \mathrm{ad}^\ell_x(y-x) \Bigg\| \\
	&\quad + \Bigg\| \sum_{\ell=1}^{r-1}
		d_\ell \sum_{\tiny \begin{array}{cc}
		(n,m)\in\mathcal{I}_\ell \\ |m|+|n|=\ell \end{array}} 
		\frac{1_{\{|n|>0\}}}{|n|+1} 
		\mathrm{ad}_x^{m_1} \mathrm{ad}_{y-x}^{n_1}
		\cdots \mathrm{ad}_x^{m_\ell} \mathrm{ad}_{y-x}^{n_\ell} (y-x)
		\Bigg\|_{\mathfrak{g}_{CM}} \\
	&\le \left(1 + \sum_{\ell=1}^{r-1}
		\sum_{\tiny \begin{array}{cc}
		(n,m)\in\mathcal{I}_\ell \\ |m|+|n|=\ell
		\end{array}} |d_\ell| C^\ell \|x\|_{\mathfrak{g}_{CM}}^\ell
		\right) \|y-x\|_{\mathfrak{g}_{CM}} +
		o\left(\|y-x\|^2_{\mathfrak{g}_{CM}}\right),
\end{align*}
where $C=C(\omega,\alpha)$ is as defined in (\ref{e.Comega}).  
Interchanging the roles of $x$ and $y$ in $g(s)$, and thus in 
this inequality, completes the proof.
\end{proof}

\begin{notation}
Let $\tau$ denote the norm topology on $G_{CM}$ and $\tau_d$ denote the
topology induced by $d_{CM}$.
\end{notation}

\begin{prop}
\label{p.length2}
For any $y\in G$ and $W\in\tau$ such that $y\in W$, 
there exists $U\in\tau_d$ such that $y\in U\subset W$. 
\end{prop}

\begin{proof} 
First we will show that, there exists $\varepsilon_0>0$ such that, 
for any $x,y\in G_{CM}$ and $\varepsilon\in(0,\varepsilon_0/2)$, 
if $d_{CM}(x,y)<\varepsilon$, then $\|x^{-1}y\|_{\mathfrak{g}_{CM}}<
2\varepsilon$.  Then we will show that the continuity
of the map $x\mapsto\|x^{-1}y\|_{\mathfrak{g}_{CM}}$ (for fixed $y$)
suffices to complete the proof.

Let $d_\ell$ be as in equation (\ref{e.dl}) and $C=C(\omega,\alpha)$ be as
in equation (\ref{e.Comega}).  Let
\[ \kappa := \sum_{\ell=1}^{r-1} |d_\ell| C^\ell, \]
and take $\varepsilon_0 := 1/2\kappa\wedge 1$.
Let $B_{\varepsilon_0} := \left\{ x\in\mathfrak{g}_{CM}:
\|x\|_{\mathfrak{g}_{CM}}\le \varepsilon_0 \right\}$.  Suppose $y\in
B_{\varepsilon_0}$, and let $g:[0,1]\rightarrow G_{CM}$ be a $C^1$-path such
that $g(0)=\mathbf{e}$ and $g(1)=y$.  Further, let $T\in[0,1]$ be the
first time that $g$ exits $B_{\varepsilon_0}$, with the convention that $T=1$ if
$g([0,1])\subset B_{\varepsilon_0}$.  Then, by equation (\ref{e.lg}), 
\begin{align*}
\ell_{CM}(g&2-) \ge \ell_{CM}(g|_{[0,T]}) \\
	&\ge \int_0^T \|g'(s)\|_{\mathfrak{g}_{CM}}
		- \sum_{\ell=1}^{r-1}
		|d_\ell| \left\| \mathrm{ad}_{g(s)}^\ell g'(s)
		\right\|_{\mathfrak{g}_{CM}}\,ds \\
	&\ge \left(1-\sum_{\ell=1}^{r-1}
		|d_\ell|C^\ell\varepsilon_0^\ell \right) 
		\int_0^T \|g'(s)\|_{\mathfrak{g}_{CM}} ds 
	\ge (1-\kappa\varepsilon_0) \|g(T)\|_{\mathfrak{g}_{CM}} 
	\ge \frac{1}{2} \|y\|_{\mathfrak{g}_{CM}}.
\end{align*} 
Taking the infimum over $g$ implies that
\[ d_{CM}(\mathbf{e},y) \ge \frac{1}{2}\|y\|_{\mathfrak{g}_{CM}},
	\qquad \text{ for all } y\in B_{\varepsilon_0}. \]
Now, if $y\notin B_{\varepsilon_0}$, then the path $g$ would have had to exit
$B_{\varepsilon_0}$ and
$\ell_{CM}(g)\ge\|g(T)\|_{\mathfrak{g}_{CM}}/2=\varepsilon_0/2$ implies 
that $d_{CM}(\mathbf{e},y)\ge\varepsilon_0/2$.  Thus,
\[ d_{CM}(\mathbf{e},y)
	\ge\frac{1}{2}\min(\varepsilon_0,\|y\|_{\mathfrak{g}_{CM}}), \qquad
	\text{ for all } y\in G_{CM}. \]
By the left invariance of $d_{CM}$, this implies that, for any $x,y\in
G_{CM}$,
\[ d_{CM}(x,y) = d_{CM}(\mathbf{e},x^{-1}y)
	\ge \frac{1}{2} \min(\varepsilon_0, \|x^{-1}y\|_{\mathfrak{g}_{CM}}). \]
So if $d_{CM}(x,y)<\varepsilon_0/2$, then $\|x^{-1}y\|_{\mathfrak{g}_{CM}}
\le 2d_{CM}(x,y) $.  

Now let $W\in\tau$ (non-empty) and fix $y\in W$.  Recall that 
Lemma \ref{l.cts} implies that
the map $x\mapsto \|x^{-1}y\|_{\mathfrak{g}_{CM}}$ is $\tau$-continuous, 
and clearly $\|x^{-1}y\|_{\mathfrak{g}_{CM}}=0$ if and only 
if $x=y$.  Thus,
\[ A_n(y) := \left\{ x : \|x^{-1}y\|_{\mathfrak{g}_{CM}}
	< \frac{1}{n} \right\} \downarrow \{y\}, \]
and there exists $N$ sufficiently large that $1/N<\varepsilon_0/2$ and 
$A_N(y)\subset W$.  Then
\[ B_N(y) : = \left\{ x: d_{CM}(x,y)< \frac{1}{2N} \right\} \in \tau_d \]
satisfies $B_N(y)\subset A_N(y)$,
since $x\in B_N(y)$ implies that $d_{CM}(x,y)<1/2N<\varepsilon_0/4$ and thus 
$\|x^{-1}y\|_{\mathfrak{g}_{CM}} \le 2 d_{CM}(x,y)< 1/N$.  
Thus, $B_N(y)\subset W$.
\end{proof} 

In particular, taking $W=\{x: \|y-x\|_{\mathfrak{g}_{CM}}<\delta\}$ for some
$\delta>0$ in this proposition, 
the proof implies that there exists $N$ such that 
$d_{CM}(x,y)<1/2N$ implies that $\|y-x\|_{\mathfrak{g}_{CM}}<\delta$.
Propositions \ref{p.length1} and \ref{p.length2} give the following corollary.

\begin{cor}
\label{c.length}
The topologies $\tau$ and $\tau_d$ are equivalent.
\end{cor}

\subsection{Ricci curvature}
\label{s.Ric}
In this section, we compute the Ricci curvature of certain
finite dimensional approximations of $G$ and show that it is bounded
below uniformly.  This result will be used in Section \ref{s.quasi}
to give $L^p$-bounds on Radon Nikodym derivatives of $\nu_t$.  
It will also be
applied in Section \ref{s.logsob} to prove a logarithmic Sobolev
inequality for $\nu_t$.  First we must define the appropriate
approximations.

Let $i:H\rightarrow W$ be the inclusion map, and $i^*:W^*\rightarrow H^*$ be
its transpose.  That is, $i^*\ell:=\ell\circ i$ for all $\ell\in W^*$.  Also,
let
\[ H_* := \{h\in H: \langle\cdot,h\rangle_H\in \mathrm{Range}(i^*)\subset H\}.
\]
That is, for $h\in H$, $h\in H_*$ if and only if $\langle\cdot,h\rangle_H\in
H^*$ extends to a continuous linear functional on $W$, which we will continue
to denote by $\langle\cdot,h\rangle_H$.  Because $H$ is a dense subspace of
$W$, $i^*$ is injective and thus has a dense range.  Since
$h\mapsto\langle\cdot,h\rangle_H$ as a map from $H$ to $H^*$ is a conjugate
linear isometric isomorphism, it follows that $H_*\ni
h\mapsto\langle\cdot,h\rangle_H\in W^*$ is a conjugate linear isomorphism
also, and so $H_*$ is a dense subspace of $H$.

Now suppose that $P:H\rightarrow H$ is a finite rank orthogonal projection
such that $PH\subset H_*$.  Let $\{k_j\}_{j=1}^m$ be an orthonormal basis for
$PH$.  Then we may extend $P$ to a (unique) continuous operator 
from $W\rightarrow H$ (still denoted by $P$) by letting
\begin{equation}
\label{e.proj} 
Pw := \sum_{j=1}^m \langle w,k_j\rangle_H k_j 
\end{equation}
for all $w\in W$.  

\begin{notation}
\label{n.proj} 
Let $\mathrm{Proj}(W)$ denote the collection of finite rank projections 
on $W$ such that $PW\subset H_*$
and $P|_{H}:H\rightarrow H$ is an orthogonal projection, that is, $P$ has 
the form
given in equation \eqref{e.proj}. Further, let $G_{P}:=PW\oplus\mathfrak{v}$ (a
subgroup of $G_{CM})$, and we equip $G_P$ with the left invariant Riemannian
metric induced from the restriction of the inner product on
$\mathfrak{g}_{CM}=H\oplus\mathfrak{v}$ to 
$\mathrm{Lie}(G_P)=PH\oplus\mathfrak{v}=:\mathfrak{g}_{CM}^P$.  Let
$\mathrm{Ric}^P$ denote the associated Ricci tensor at the identity in $G_P$.
\end{notation}

\begin{prop}
\label{p.Ric2}
For $X=(A,a)\in\mathfrak{g}_{CM}^P$,
\[ \langle \mathrm{Ric}^P X,X\rangle_{\mathfrak{g}^P_{CM}} 
	= \frac{1}{4} \|\langle a,[\cdot,\cdot]\rangle\|_{(\mathfrak{g}_{CM}^P)^*
		\otimes (\mathfrak{g}_{CM}^P)^*}^2 
	-\frac{1}{2} \|[\cdot,X]\|_{(\mathfrak{g}_{CM}^P)^*
		\otimes \mathfrak{v}}^2,\]
where $(\mathfrak{g}_{CM}^P)^*=(PH)^*\otimes\mathfrak{v}^*$.
\end{prop}

\begin{proof}
For $\mathfrak{g}$ any nilpotent Lie algebra with orthonormal basis
$\Gamma$, 
\begin{equation}
\label{e.Ric0}
\langle \mathrm{Ric}\, X,X\rangle = \frac{1}{4}\sum_{Y\in\Gamma}
\|\mathrm{ad}^*_Y X\|^2 - \frac{1}{2}\sum_{Y\in\Gamma} \|\mathrm{ad}_Y X\|^2,
\end{equation}
for all $X\in\mathfrak{g}$;
see for example Theorem 7.30 and Corollary 7.33 of \cite{Besse87}. 

So let $\Gamma_m:=\{h_i\}_{i=1}^{m+N} = 
\{(k_i,0)\}_{i=1}^m\cup\{(0,e_j)\}_{j=1}^N$ be an
orthonormal basis of $\mathfrak{g}_{CM}^P=PH\oplus\mathfrak{v}$, where
$\{k_i\}_{i=1}^m$ and $\{e_j\}_{j=1}^N$ are orthonormal bases of $PH$ and
$\mathfrak{v}$, respectively.  Then, for $Y\in\mathfrak{g}_{CM}^P$,
\[ \mathrm{ad}_Y^*X 
	= \sum_{h_i\in\Gamma_m} \langle \mathrm{ad}_Y^* X, h_i\rangle
		_{\mathfrak{g}_{CM}} h_i 
	= \sum_{h_i\in\Gamma_m} \langle X, \mathrm{ad}_Y h_i\rangle 
		_{\mathfrak{g}_{CM}} h_i. \]
Thus,
\[ \sum_{h_i\in\Gamma_m} \|\mathrm{ad}_{h_i}^*X\|_{\mathfrak{g}_{CM}}^2
	= \sum_{h_i\in\Gamma_m}\sum_{h_j\in\Gamma_m}  
		\langle X, \mathrm{ad}_{h_i} h_j \rangle_{\mathfrak{g}_{CM}}^2 
	= \sum_{h_i,h_j\in\Gamma_m} \langle X, [h_i,h_j] \rangle
		_{\mathfrak{g}_{CM}} ^2. \]
Plugging this into (\ref{e.Ric0}) gives
\begin{align*} 
\langle \mathrm{Ric}^P X,X\rangle_{\mathfrak{g}^P_{CM}} 
	&= \frac{1}{4} \sum_{h_i,h_j\in\Gamma_m} \langle X, [h_i,h_j] \rangle
		_{\mathfrak{g}_{CM}}^2
		- \frac{1}{2} \sum_{h_i\in\Gamma_m} \|[h_i,X]\|
		_{\mathfrak{g}_{CM}}^2 \\
	&= \frac{1}{4} \sum_{h_i,h_j\in\Gamma_m} \langle a, [h_i,h_j] \rangle
		_\mathfrak{v}^2
		- \frac{1}{2} \sum_{h_i\in\Gamma_m} \|[h_i,X]\|_\mathfrak{v}^2.
\end{align*}
\end{proof}

\begin{cor}
\label{c.Ric3}
Let
\[ K := -\frac{1}{2} \sup \left\{ 
		\|[\cdot,X]\|^2_{\mathfrak{g}_{CM}^*\otimes\mathfrak{v}}   
		: \,\|X\|_{\mathfrak{g}_{CM}}=1 \right\}. \]
Then $K>-\infty$ and $K$ is the largest constant such that
\[ \langle \mathrm{Ric}^P X,X\rangle_{\mathfrak{g}_{CM}^P} \ge 
	K \|X\|^2_{\mathfrak{g}_{CM}^P}, 
	\quad \text{ for all } X \in\mathfrak{g}_{CM}^P, \]
holds uniformly for all $P\in\mathrm{Proj}(W)$.
\end{cor}
\begin{proof}
The first assertion is simple, since
\[ K\ge - \frac{1}{2} \|[\cdot,\cdot]\|_2^2 > -\infty, \]
by Corollary \ref{c.normest}.  Now, for $P\in\mathrm{Proj}(W)$ as in 
Notation \ref{n.proj}, Proposition \ref{p.Ric2} implies that
\[
\langle \mathrm{Ric}^P X, X\rangle_{\mathfrak{g}_{CM}^P} 
	\ge - \frac{1}{2} \|[\cdot,X]\|^2
	_{(\mathfrak{g}_{CM}^P)^*\otimes\mathfrak{v}}.
\]
Thus,
\begin{align} 
\notag
\frac{\langle \mathrm{Ric}^P X,X \rangle_{\mathfrak{g}_{CM}^P}}
		{\|X\|^2_{\mathfrak{g}_{CM}^P}}
	&\ge - \frac{1}{2} \frac{\|[\cdot,X]\|^2_{(\mathfrak{g}^P_{CM})^*
		\otimes\mathfrak{v}}}{\|X\|^2_{\mathfrak{g}_{CM}^P}} \\
	&\label{e.pah}
	\ge - \frac{1}{2}\sup \left\{ 
		\|[\cdot,X]\|^2_{(\mathfrak{g}_{CM}^P)^*\otimes\mathfrak{v}}
 		:\quad \|X\|_{\mathfrak{g}_{CM}^P}=1 \right\}
	=: K_P.
\end{align}
Noting that the infimum of $K_P$ over all $P\in\mathrm{Proj}(W)$ is
$K$ completes the proof.
\end{proof}

\begin{remark}
Of course, one can compute the Ricci curvature for
$\mathfrak{g}=W\oplus\mathfrak{v}$ just as in Proposition \ref{p.Ric2}.
Choose an orthonormal basis $\Gamma=\{h_i\}_{i=1}^\infty =
\{(k_i,0)\}_{i=1}^\infty\cup\{(0,e_j)\}_{j=1}^N$
of $\mathfrak{g}_{CM} = H\oplus \mathfrak{v}$, 
where
$\{k_i\}_{i=1}^\infty$ is an orthonormal basis of $H$, and $\{e_j\}_{j=1}^N$
is an orthonormal basis of $\mathfrak{v}$.  Then, for all 
$X=(A,a)\in\mathfrak{g}_{CM}$,
\begin{align*} 
\langle \mathrm{Ric}\, X,X\rangle_{\mathfrak{g}_{CM}} 
	&= \frac{1}{4} \sum_{i,j=1}^\infty \langle a, [h_i,h_j] \rangle
		_\mathfrak{v}^2
		- \frac{1}{2} \sum_{i=1}^\infty \|[h_i,X]\|_\mathfrak{v}^2 \\
	&= \frac{1}{4} \|\langle
a,[\cdot,\cdot]\rangle\|_{\mathfrak{g}_{CM}^*
		\otimes \mathfrak{g}_{CM}^*}^2 
	-\frac{1}{2} \|[\cdot,X]\|_{\mathfrak{g}_{CM}^*
		\otimes \mathfrak{v}}^2
	\ge K\|X\|_{\mathfrak{g}_{CM}}^2.
\end{align*}
\end{remark}

\section{Brownian motion}
\label{s.BM}
Suppose that $B_t$ is a smooth curve in $\mathfrak{g}_{CM}$ with
$B_0=0$, and consider the differential equation
\[ \dot{g}_t = L_{g_t*}\dot{B}_t, \quad \text{ with } g_0=\mathbf{e}. \]  
The solution $g_t$ may be written as follows (see
\cite{Strichartz87}):  For $t>0$, let $\Delta_n(t)$ denote the simplex in
$\mathbb{R}^n$ given by
\[ \{s=(s_1,\cdots,s_n)\in\mathbb{R}^n: 0<s_1<s_2<\cdots<s_n<t\}. \]
Let $\mathcal{S}_n$ denote the permutation group on $(1,\cdots,n)$, and, for each
$\sigma\in \mathcal{S}_n$, let $e(\sigma)$ denote the number of ``errors'' in the
ordering $(\sigma(1),\sigma(2),\cdots,\sigma(n))$, that is, $e(\sigma)=\#
\{j<n: \sigma(j)>\sigma(j+1)\}$.  Then
\begin{multline}
\label{e.ode}
g_t = \sum_{n=1}^r \sum_{\sigma\in \mathcal{S}_n} 
	\left( (-1)^{e(\sigma)}\bigg/ n^2 
	\begin{bmatrix} n-1 \\ e(\sigma) \end{bmatrix}\right) \times \\
	\int_{\Delta_n(t)} [ \cdots[\dot{B}_{s_{\sigma(1)}},
	\dot{B}_{s_{\sigma(2)}}],\ldots,] \dot{B}_{s_{\sigma(n)}}]\, ds.
\end{multline}
For $n\in\{1,\cdots,r\}$ and $\sigma\in\mathcal{S}_n$,
let $F_n^\sigma:\mathfrak{g}_{CM}^{\otimes n}
\rightarrow\mathfrak{v}$ be the linear map given by
\begin{equation}
\label{e.Fn}
F_n^\sigma(k_1\otimes\cdots\otimes k_n) 
	:= [ [\cdots[k_{\sigma(1)},k_{\sigma(2)}],\cdots], 
		k_{\sigma(n)}].
\end{equation}
Recall that $F_n^\sigma$ is Hilbert-Schmidt by Corollary \ref{c.normest}.
Then we may write
\begin{equation}
\label{e.dotg} 
g_t = \sum_{n=1}^{r-1} \sum_{\sigma\in\mathcal{S}_n} 
	c^\sigma_n F^\sigma_n \left(
	\int_{\Delta_n(t)} \dot{B}_{s_1}\otimes\cdots\otimes
	\dot{B}_{s_n} \,ds\right) .
\end{equation}
Using this as our motivation, we first explore stochastic integral
analogues of equation (\ref{e.dotg}) where the smooth curve $B$ is replaced by
Brownian motion on $\mathfrak{g}$.

\subsection{Multiple It\^o integrals}

Let $\langle \cdot, \cdot\rangle_{\mathfrak{g}_{CM}^{\otimes n}}$ denote the inner
product on $\mathfrak{g}_{CM}^{\otimes n}$ arising from the inner product on
$\mathfrak{g}_{CM}$.  Also, let $\{k_i\}_{i=1}^\infty\subset H_*$ be an
orthonormal basis of $H$, and define $P_m\in\mathrm{Proj}(W)$ by
\begin{equation}
\label{e.pm}
P_m(w) = \sum_{i=1}^m \langle w,k_i\rangle_H k_i, 
	\qquad \text{ for all } w\in W,
\end{equation}
as in equation (\ref{e.proj}), and define
\begin{equation}
\label{e.pim}
\pi_m(w,x) := \pi_{P_m}(w,x) := (P_m(w),x) \in G_{P_m}.
\end{equation}
Of course, $\mathrm{dim}(G_{P_m})=m+N$, but
in a mild abuse of notation, we will use $\{h_i\}_{i=1}^m$ to denote 
an orthonormal basis of $G_{P_m}$, rather than
the more cumbersome $\{h_i\}_{i=1}^{m+N}=\{(k_i,0)\}_{i=1}^m\cup
\{(0,e_i)\}_{i=1}^N$, where $\{e_i\}_{i=1}^N$ is an orthonormal basis of
$\mathfrak{v}$.  

Let $\{B_t\}_{t\ge0}=\{(\beta_t,\beta^\mathfrak{v}_t)\}_{t\ge0}$ 
be a Brownian motion on $\mathfrak{g}=W\oplus\mathfrak{v}$ 
with variance determined by
\[
\mathbb{E}\left[\langle B_s,h\rangle
		_{\mathfrak{g}_{CM}} 
		\langle B_t,k\rangle
		_{\mathfrak{g}_{CM}}\right] 
	= \langle h,k \rangle_{\mathfrak{g}_{CM}} \min(s,t),
\]
for all $s,t\ge0$ and $h=(A,a)$ and $k=(C,c)$, such that
$A,C\in H_*$ and $a,c\in\mathfrak{v}$.
Then $\pi_m B= (P_m\beta,\beta^\mathfrak{v})$ is a Brownian motion on
$\mathfrak{g}^{P_m}=P_mW\oplus\mathfrak{v}\subset\mathfrak{g}_{CM}$.

\begin{prop}
\label{p.int}
For $\xi\in L^2(\Delta_n(t),\mathfrak{g}_{CM}^{\otimes n})$ a continuous 
mapping, let
\[ 
J_n^m(\xi)_t := \int_{\Delta_n(t)} \langle \xi(s), d\pi_m B_{s_1}
	\otimes\cdots\otimes d\pi_m B_{s_n}
	\rangle_{\mathfrak{g}_{CM}^{\otimes n}}. 
\]
Then $\{J_n^m(\xi)_t\}_{t\ge0}$ is a continuous $L^2$-martingale such that,
for all $m$,
\[ \mathbb{E}|J_n^m(\xi)_t|^2\le
	\|\xi\|^2_{L^2(\Delta_n(t),\mathfrak{g}_{CM}^{\otimes n})}, \] 
and there exists a
continuous $L^2$-martingale $\{J_n(\xi)_t\}_{t\ge0}$ such that
\begin{equation} 
\label{e.Jnm}
\lim_{m\rightarrow\infty} \mathbb{E}\left[ \sup_{\tau\le t} 
	|J_n^m(\xi)_\tau-J_n(\xi)_\tau|^2 \right] = 0,
\end{equation}
for all $t<\infty$.  In particular,
\begin{equation}
\label{e.mart}
J_n(\xi)_t 
	:= \int_{\Delta_n(t)} \langle \xi(s), dB_{s_1}\otimes\cdots\otimes dB_{s_n}
	\rangle_{\mathfrak{g}_{CM}^{\otimes n}},
\end{equation}
and $J_n(\xi)_t$ is well-defined independent of 
the choice of orthonormal basis $\{h_i\}_{i=1}^\infty$ in (\ref{e.pm}).
\end{prop}

\begin{proof} 
Note first that, 
\[ J_n^m(\xi)_t = \sum_{i_1,\ldots,i_n=1}^m 
	\int_{\Delta_n(t)} \langle\xi(s),h_{i_1}\otimes\cdots\otimes
	h_{i_n}\rangle_{\mathfrak{g}_{CM}^{\otimes n}}
	dB_{s_1}^{i_1}\cdots dB_{s_n}^{i_n} 
\]
where $\{B^i\}_{i=1}^m$ are independent real
valued Brownian motions.  Let $\xi_{i_1,\ldots,i_n} :=
\langle\xi,h_{i_1}\otimes\cdots\otimes h_{i_n}\rangle$.  Then
\[ |\xi_{i_1,\ldots,i_n}(s) |^2\le \|\xi(s)\|^2
	_{\mathfrak{g}_{CM}^{\otimes n}} \]
and $\xi_{i_1,\ldots,i_n} \in L^2(\Delta_n(t))$.  Thus, $J_n^m(\xi)_t$ is 
defined as a
(finite dimensional) vector-valued multiple Wiener-It\^o integral, see for
example \cite{Ito51,Shigekawa04}.

Now note that
\begin{align*} 
dJ&_n^m(\xi)_t 
	= \int_{\Delta_{n-1}(t)} \langle \xi(s_1,\ldots,s_{n-1},t),
		d\pi_mB_{s_1}\otimes\cdots\otimes d\pi_mB_{s_{n-1}} \otimes d\pi_mB_t
		\rangle_{\mathfrak{g}_{CM}^{\otimes n}} \\
	&= \sum_{i=1}^m \int_{\Delta_{n-1}(t)} \langle \xi(s_1,\ldots,s_{n-1},t),
		d\pi_mB_{s_1}\otimes\cdots\otimes d\pi_mB_{s_{n-1}} \otimes h_i
		\rangle_{\mathfrak{g}_{CM}^{\otimes n}} dB_t^i.
\end{align*}
Thus, the quadratic variation $\langle J_n^m(\xi) \rangle_t$ is given by
\[ \sum_{i=1}^m \int_0^t 
	\bigg|\int_{\Delta_{n-1}(\tau)}
	\langle \xi(s_1,\ldots,s_{n-1},\tau),
	d\pi_mB_{s_1}\otimes\cdots\otimes
	d\pi_mB_{s_{n-1}}\otimes h_i \rangle
	_{\mathfrak{g}_{CM}^{\otimes n}}\bigg|^2 d\tau,
\]
and
\begin{align*}
\mathbb{E}&|J_n^m(\xi)_t|^2 
	= \mathbb{E} \langle J_n^m(\xi) \rangle_t \\
	&= \sum_{i_1=1}^m \int_0^t 
		\mathbb{E}\bigg[\sum_{i_2=1}^m \int_0^{\tau_1}
		\bigg|\int_{\Delta_{n-2}(\tau_2)}
		\langle \xi(s_1,\ldots,s_{n-2},\tau_2,\tau_1),
		d\pi_mB_{s_1}\otimes\cdots \\
	&\qquad\qquad\qquad\qquad\qquad\qquad\qquad\qquad 
		\cdots \otimes d\pi_mB_{s_{n-2}}
		\otimes h_{i_2}\otimes h_{i_1}\rangle
		_{\mathfrak{g}_{CM}^{\otimes n}}\bigg|^2d\tau_2 \bigg]d\tau_1.
\end{align*}
Iterating this procedure $n$ times gives
\begin{align}
\label{e.gab}
\mathbb{E}|J_n^m(\xi)_t|^2 
	&= \sum_{i_1,\ldots,i_n=1}^m \int_{\Delta_n(t)} 
		\left|\langle \xi(\tau_1,\cdots,\tau_n), 
		h_{i_1}\otimes\cdots\otimes h_{i_n} \rangle
		_{\mathfrak{g}_{CM}^{\otimes n}}\right|^2 d\tau_1\cdots d\tau_n \\
	&\notag
	= \int_{\Delta_n(t)} \|\pi_m^{\otimes n} \xi(s) \|
		_{\mathfrak{g}_{CM}^{\otimes n}}^2
	\le \|\xi \|_{L^2(\Delta_n(t),\mathfrak{g}_{CM}^{\otimes n})}^2,
\end{align}
and thus, for each $n$, $J_n^m(\xi)_t$ is bounded uniformly in $L^2$
independent of $m$.  

A similar argument shows that the sequence 
$\left\{J_n^m(\xi)_t\right\}_{m=1}^\infty$ is Cauchy in $L^2$.  
For $m\le\ell$, consider
\begin{multline}
\label{e.gub} 
J_n^\ell(\xi)_t - J_n^m(\xi)_t = 
	\sum_{j=1}^{n} \int_{\Delta_n(t)} \langle \xi(s) , 
	d\pi_\ell B_{s_1} \otimes \cdots \\ \cdots\otimes d\pi_\ell B_{s_{j-1}}
	\otimes d(\pi_\ell-\pi_m)B_{s_j} \otimes d\pi_m B_{s_{j+1}}
	\otimes\cdots\otimes d\pi_m B_{s_n} \rangle. 
\end{multline}
Thus, applying Cauchy-Schwarz and computing as in equation (\ref{e.gab}),
\begin{multline}
\label{e.tod}
\mathbb{E}\left|J_n^\ell(\xi)_t - J_n^m(\xi)_t \right|^2 \\
	\le n \sum_{j=1}^n \sum_{i_1,\ldots,i_{j-1}=1}^\ell \sum_{i_j=m+1}^\ell
		\sum_{i_{j+1},\ldots,i_n=1}^m \int_{\Delta_n(t)} 
		\left|\langle \xi(s), h_{i_1}\otimes\cdots\otimes h_{i_n} \rangle
		_{\mathfrak{g}_{CM}^{\otimes n}}\right|^2 ds 
	\rightarrow 0,
\end{multline}
as $\ell,m\rightarrow \infty$, since 
\begin{align*} 
\|\xi\|_{L^2(\Delta_n(t),\mathfrak{g}_{CM}^{\otimes n})}^2
	&= \int_{\Delta_n(t)} \|\xi(s)\|_{\mathfrak{g}_{CM}^{\otimes n}}^2\,ds \\
	&= \int_{\Delta_n(t)} \sum_{i_1,\ldots,i_n=1}^\infty
		\left|\langle \xi(s), h_{i_1}\otimes\cdots\otimes h_{i_n}\rangle
		_{\mathfrak{g}_{CM}^{\otimes n}}\right|^2\,ds 
	<\infty. 
\end{align*}
Since the space of continuous $L^2$-martingales is complete in the
norm $M\mapsto \mathbb{E}|M_t|^2$, there exists a continuous martingale
$\{X_t\}_{t\ge0}$ such that
\begin{equation}
\label{e.tad} 
\lim_{m\rightarrow\infty} \mathbb{E}|J_n^m(\xi)_t-X_t|^2 = 0. 
\end{equation}

To see that $X_t$ is independent of basis, suppose now that
$\{h'_j\}_{j=1}^\infty\subset H_*$ is another orthonormal
basis for $H$ and $P_m':W\rightarrow H_*$ and $\pi_m':G\rightarrow G_{P_m'}$ 
are the corresponding orthogonal projections, that is,
\[ P_m'w := \sum_{i=1}^m \langle w,h_i'\rangle_W h_i', \]
and $\pi_m'(w,x) = (P_m'w,x)$.  Let 
\[ 
J_n^{m'}(\xi)_t = \int_{\Delta_n(t)} \langle \xi(s), d\pi_m' B_{s_1}
	\otimes\cdots\otimes d\pi_m' B_{s_n}
	\rangle_{\mathfrak{g}_{CM}^{\otimes n}}.
\]
Then, using equation (\ref{e.gub}) with $\pi_\ell$ replaced by $\pi_m'$,
applying Cauchy-Schwarz, and again computing as in (\ref{e.gab}), gives
\begin{multline*}
\mathbb{E}\left|J_n^m(\xi)_t - J_n^{m'}(\xi)_t \right|^2 
	\le n  \sum_{j=1}^n \int_{\Delta_n(t)} \sum_{i_1,\ldots,i_n=1}^\infty
		\bigg|\langle \xi(s), \pi_mh_{i_1}\otimes\cdots \\
		\cdots\otimes \pi_m h_{i_{j-1}} \otimes (\pi_m-\pi_m') h_{i_j} 
		\otimes \pi_m' h_{i_{j+1}}
		\otimes\cdots\otimes \pi_m' h_{i_n} \rangle
		_{\mathfrak{g}_{CM}^{\otimes n}}\bigg|^2 ds.
\end{multline*}
Writing $\pi_m-\pi_m'=(\pi_m-I) + (I-\pi_m')$, and 
considering terms for each fixed $j$, we have 
\begin{align*}
\int_{\Delta_n(t)}& \sum_{i_1,\ldots,i_n=1}^\infty
		\bigg|\langle \xi(s), \pi_m h_{i_1}\otimes
		\cdots\otimes \pi_m h_{i_{j-1}} \otimes (\pi_m-I) h_{i_j} \\
	&\qquad\qquad\qquad\qquad\qquad\qquad\qquad\qquad
		\otimes \pi_m' h_{i_{j+1}}
		\otimes\cdots\otimes \pi_m' h_{i_n} \rangle
		_{\mathfrak{g}_{CM}^{\otimes n}}\bigg|^2 ds \\
	&= \int_{\Delta_n(t)} \sum_{i_1,\ldots,i_j=1}^m \sum_{i_j=m+1}^\infty
		\sum_{i_{j+1},\ldots,i_n=1}^\infty
		\bigg|\langle \xi(s), h_{i_1}\otimes
		\cdots\otimes h_{i_{j-1}} \otimes h_{i_j} \\
	&\qquad\qquad\qquad\qquad\qquad\qquad\qquad\qquad\qquad
		\otimes \pi_m' h_{i_{j+1}}
		\otimes\cdots\otimes \pi_m' h_{i_n} \rangle
		_{\mathfrak{g}_{CM}^{\otimes n}}\bigg|^2 ds \\
	&\le \int_{\Delta_n(t)} \sum_{i_1,\ldots,i_j=1}^m \sum_{i_j=m+1}^\infty
		\sum_{i_{j+1},\ldots,i_n=1}^\infty
		\left|\langle \xi(s), h_{i_1}\otimes
		\cdots\otimes h_{i_n} \rangle
		_{\mathfrak{g}_{CM}^{\otimes n}}\right|^2 ds
	\rightarrow 0,
\end{align*}
as $m\rightarrow\infty$.  Similarly,
\begin{align*}
&\int_{\Delta_n(t)} \sum_{i_1,\ldots,i_n=1}^\infty
	\bigg|\langle \xi(s), \pi_m h_{i_1}\otimes
		\cdots\otimes \pi_m h_{i_{j-1}} \\
	&\qquad\qquad\qquad\qquad\qquad\qquad\otimes (\pi_m'-I) h_{i_j} 
		\otimes \pi_m' h_{i_{j+1}}
		\otimes\cdots\otimes \pi_m' h_{i_n} \rangle
		_{\mathfrak{g}_{CM}^{\otimes n}}\bigg|^2 ds \\
	&\quad= \int_{\Delta_n(t)} \sum_{i_1,\ldots,i_n=1}^\infty	
		\bigg|\langle \xi(s), \pi_mh_{i_1}'\otimes
		\cdots\otimes \pi_mh_{i_{j-1}}' \\
	&\qquad\qquad\qquad\qquad\qquad\qquad \otimes (\pi_m'-I) h_{i_j}' 
		\otimes \pi_m' h_{i_{j+1}}'
		\otimes\cdots\otimes \pi_m' h_{i_n}' \rangle
		_{\mathfrak{g}_{CM}^{\otimes n}} \bigg|^2 ds \\
	&\quad= \int_{\Delta_n(t)} \sum_{i_1,\ldots,i_{j-1}=1}^\infty 
		\sum_{i_j=m+1}^\infty \sum_{i_{j+1},\ldots,i_n=1}^m
		\bigg|\langle \xi(s), \pi_mh_{i_1}'\otimes
		\cdots\otimes \pi_mh_{i_{j-1}}' \\
	&\qquad\qquad\qquad\qquad\qquad\qquad \otimes h_{i_j}' 
		\otimes h_{i_{j+1}}' \otimes\cdots\otimes h_{i_n}' \rangle
		_{\mathfrak{g}_{CM}^{\otimes n}} \bigg|^2 ds \\
	&\quad\le \int_{\Delta_n(t)} \sum_{i_1,\ldots,i_{j-1}=1}^\infty 
		\sum_{i_j=m+1}^\infty \sum_{i_{j+1},\ldots,i_n=1}^m
		\left|\langle \xi(s), h_{i_1}'\otimes\cdots\otimes h_{i_n}' \rangle
		_{\mathfrak{g}_{CM}^{\otimes n}} \right|^2 ds 
	\rightarrow 0,
\end{align*}
as $m\rightarrow\infty$.  Thus,
\[ \lim_{m\rightarrow\infty} \mathbb{E}\left|J_n^m(\xi)_t 
	- J_n^{m'}(\xi)_t \right|^2 = 0, \]
and $X_t$ is independent of the choice of orthonormal basis.  
In particular, replacing $J_n^\ell(\xi)_t$ in (\ref{e.tod}) by
$J_n(\xi)_t$ as given in equation (\ref{e.mart}), and taking the
limit as $m\rightarrow\infty$,  shows that $X_t=J_n(\xi)_t$ 
satisfies (\ref{e.tad}).  Combining this with Doob's
maximal inequality proves equation (\ref{e.Jnm}).

\end{proof}

A simple linearity argument extends the map $J_n$ to functions taking
values in $(\mathfrak{g}_{CM}^*)^{\otimes n}\otimes\mathfrak{v}$.

\begin{cor}
\label{c.int}
Let $F\in L^2(\Delta_n(t),(\mathfrak{g}_{CM}^*)^{\otimes
n}\otimes\mathfrak{v})$ be a continuous map.  That is, 
$F:\Delta_n(t)\times\mathfrak{g}_{CM}^{\otimes n}\rightarrow
\mathfrak{v}$ is a map continuous in $s$ and linear on
$\mathfrak{g}_{CM}^{\otimes n}$ such that 
\[ \int_{\Delta_n(t)} \|F(s)\|_2^2 \,ds
	= \int_{\Delta_n(t)} \sum_{j_1,\ldots,j_n=1}^\infty 
		\|F(s)(h_{j_1}\otimes\cdots\otimes h_{j_n})\|_\mathfrak{v}^2\,ds
	<\infty. \]
Then
\[ J_n^m(F)_t := \int_{\Delta_n(t)} F(s)(d\pi_m B_{s_1}
	\otimes\cdots\otimes d\pi_m B_{s_n}) 
\]
is a continuous $L^2$-martingale, and there exists a continuous 
$\mathfrak{v}$-valued $L^2$-martingale $\{J_n(F)_t\}_{t\ge0}$ 
such that
\[
\lim_{m\rightarrow\infty} \mathbb{E}\left[ \sup_{\tau\le t} 
	\|J_n^m(\xi)_\tau-J_n(\xi)_\tau\|_\mathfrak{v}^2 \right] = 0,
\]
for all $t<\infty$.  The martingale $J_n(\xi)_t$ is well-defined independent of 
the choice of orthonormal basis $\{h_i\}_{i=1}^\infty$ in (\ref{e.pm}), 
and will be denoted by
\[
J_n(F)_t := \int_{\Delta_n(t)} F(s)(dB_{s_1}\otimes\cdots\otimes dB_{s_n}).
\]
\end{cor}
\begin{proof}
Let $\{e_j\}_{j=1}^N$ be an orthonormal basis of $\mathfrak{v}$.  Then for any
$k_1,\ldots,k_n\in\mathfrak{g}_{CM}$,
\[ F(s)(k_1\otimes\cdots\otimes k_n) 
	= \sum_{j=1}^N \langle F(s) (k_1\otimes\cdots\otimes k_n), 
		e_j \rangle e_j. \] 
Since $\langle F(s)(\cdot), e_j \rangle$ is linear on
$\mathfrak{g}_{CM}^{\otimes n}$, for each $s$ there exists
$\xi_j(s)\in\mathfrak{g}_{CM}^{\otimes n}$ such that
\begin{equation}
\label{e.xij} 
\langle \xi_j(s), k_1\otimes\cdots\otimes k_n \rangle 
	= \langle F(s)(k_1\otimes\cdots\otimes k_n), e_j \rangle. 
\end{equation}
If $\xi_j:\Delta_n(t)\rightarrow\mathfrak{g}_{CM}^{\otimes n}$ is defined
by equation (\ref{e.xij}), then 
\[ \|\xi_j\|_{L^2(\Delta_n(t),\mathfrak{g}_{CM}^{\otimes n})} 
	\le \int_{\Delta_n(t)} \|F(s)\|_2^2\, ds < \infty. \]
Thus,
\[ J_n(F)_t = \sum_{j=1}^N \int_{\Delta_n(t)} \langle \xi_j(s),
		dB_{s_1}\otimes\cdots\otimes dB_{s_n} \rangle e_j
	= \sum_{j=1}^N J_n(\xi_j)_t e_j, \]
is well-defined, and, for each $j$, $J_n(\xi_j)$ is a martingale as defined in 
Proposition \ref{p.int}.

\end{proof}

\subsection{Brownian motion and finite dimensional approximations}

Again let $B_t$ denote Brownian motion on $\mathfrak{g}$.
By equation (\ref{e.ode}), the solution to the Stratonovich
stochastic differential equation
\[ \delta g_t = L_{g_t*} \delta B_t, \quad \text{ with } g_0=\mathbf{e}, \]
should be given by
\begin{align*}
g_t = \sum_{n=1}^{r-1} \sum_{\sigma\in\mathcal{S}_n} c^\sigma_n \int_{\Delta_n(t)} 
	[ [\cdots[\delta B_{s_{\sigma(1)}},\delta B_{s_{\sigma(2)}}],\cdots], 
	\delta B_{s_{\sigma(n)}}],
\end{align*}
for coefficients $c_n^\sigma$ determined by equation (\ref{e.ode}).

To see that this process is well-defined, consider the following.
Let $\{M_n(t)\}_{t\ge0}$ denote the process in $\mathfrak{g}^{\otimes n}$ 
defined by
\[ 
M_n(t) := \int_{\Delta_n(t)} \delta B_{s_1}\otimes\cdots\otimes \delta
	B_{s_n}.
\]
By repeatedly applying the definition of the Stratonovich integral, 
the iterated Stratonovich integral $M_n(t)$ 
may be realized as a linear combination of iterated It\^o integrals:
\[ M_n(t) = \sum_{m=\lceil n/2\rceil}^n \frac{1}{2^{n-m}}
		\sum_{\alpha\in\mathcal{J}_n^m} I^n_t(\alpha), \]
where
\[ \mathcal{J}_n^m := \left\{(\alpha_1,\ldots,\alpha_m)\in\{1,2\}^m :
	\sum_{i=1}^m \alpha_i = n \right\}, \]
and, for $\alpha\in\mathcal{J}_n^m$, $I_t^n(\alpha)$ is the iterated It\^o integral
\[ I_t^n(\alpha) = \int_{\Delta_m(t)} dX^1_{s_1}\otimes\cdots\otimes 
	dX^m_{s_m} \]
with
\[ dX_s^i = \left\{ \begin{array}{cl} dB_s & \text{if } \alpha_i=1 \\
	\sum_{j=1}^\infty h_j \otimes h_j \, ds & \text{if } \alpha_i=2
	\end{array} \right. ; \]
compare with Proposition 1 of \cite{BenArous89}.

As in equation (\ref{e.Fn}), letting 
\[ 
F_n^\sigma(k_1\otimes\cdots\otimes k_n) 
	:= [ [\cdots[k_{\sigma(1)},k_{\sigma(2)}],\cdots], 
		k_{\sigma(n)}],
\]
we may write
\begin{align*}
g_t &= \sum_{n=1}^{r-1} \sum_{\sigma\in\mathcal{S}_n} 
	c^\sigma_n F^\sigma_n (M_n(t)) \\
	&= \sum_{n=1}^{r-1} \sum_{\sigma\in\mathcal{S}_n}
	\sum_{m=\lceil n/2\rceil}^n \frac{c^\sigma_n }{2^{n-m}}
		\sum_{\alpha\in\mathcal{J}_n^m} F^\sigma_n (I^n_t(\alpha)), 
\end{align*}
presuming the integrals $F_n^\sigma(I_t^n(\alpha))$ are defined.

For each $\alpha$, let $p_\alpha=\#\{i:\alpha_i=1\}$ and $q_\alpha=\#\{i:
\alpha_i=2\}$ (so that $p_\alpha+q_\alpha=m$ when $\alpha\in\mathcal{J}_n^m$),
and let 
\[ \mathcal{J}_n := \bigcup_{m=\lceil n/2\rceil}^n \mathcal{J}_n^m. \]
Then, for each $\sigma\in\mathcal{S}_n$ and $\alpha\in\mathcal{J}_n$, 
\[ F_n^\sigma(I_t^n(\alpha))
	= \int_{\Delta_{p_\alpha}(t)} f_\alpha(s,t) \hat{F}_n^{\sigma,\alpha}
		(dB_{s_1}\otimes\cdots\otimes dB_{s_{p_\alpha}}),
\]
where $\hat{F}_n^{\sigma,\alpha}$ and $f_\alpha$ are defined as follows.
$\hat{F}_n^{\sigma,\alpha}:\mathfrak{g}^{\otimes p_\alpha}
\rightarrow\mathfrak{g}$ is defined by
\begin{multline}
\label{e.hatF} 
\hat{F}_n^{\sigma,\alpha}(k_1\otimes\cdots\otimes k_{p_\alpha}) \\
	:= \sum_{j_1,\ldots,j_{q_\alpha}=1}^\infty 
		F_n^{\sigma'}(k_1\otimes\cdots\otimes k_{p_\alpha}
		\otimes h_{j_1}\otimes h_{j_1}
		\otimes\cdots\otimes h_{j_{q_\alpha}}\otimes h_{j_{q_\alpha}}), 
\end{multline}
for $\{h_j\}_{j=1}^\infty$ an orthonormal basis of $\mathfrak{g}_{CM}$ and
$\sigma'=\sigma'(\alpha)\in\mathcal{S}_n$ given by 
$\sigma'=\sigma\circ\tau^{-1}$, for any $\tau\in\mathcal{S}_n$ such that
\begin{multline*} 
\tau(dX^1_{s_1}\otimes\cdots\otimes dX^m_{s_m}) \\
	= \sum_{j_1,\cdots,j_{q_\alpha}=1}^\infty dB_{s_1}\otimes\cdots
	\otimes dB_{s_{p_\alpha}}\otimes h_{j_1}\otimes h_{j_1}\otimes\cdots
	\otimes h_{j_{q_\alpha}}\otimes h_{j_{q_\alpha}} ds_1\cdots
	ds_{q_\alpha}.
\end{multline*}
The function $f_\alpha$ is a polynomial of order
$q_\alpha$ in $s=(s_1,\ldots,s_{p_\alpha})$ and $t$.  
Thus, $f_\alpha$ may be written as
\begin{equation}
\label{e.falpha} 
f_\alpha(s,t) = \sum_{a=0}^{q_\alpha} b_\alpha^a t^a 
	\tilde{f}_{\alpha,a}(s), 
\end{equation}
for some coefficients $b_\alpha^a\in\mathbb{R}$ and polynomials 
$\tilde{f}_{\alpha,a}$ of degree $q_\alpha-a$ in $s$.  
If $\hat{F}_n^{\sigma,\alpha}$ is
Hilbert-Schmidt on $\mathfrak{g}_{CM}^{\otimes p_\alpha}$, then
\[ \int_{\Delta_{p_\alpha}(t)} \left\|\tilde{f}_{\alpha,a}(s)
	\hat{F}_n^{\sigma,\alpha} \right\|_2^2\, ds	
	= \left\|\tilde{f}_{\alpha,a}\right\|_{L^2(\Delta_{p_\alpha}(t))}
		\left\|\hat{F}_n^{\sigma,\alpha} \right\|_2^2
	<\infty,
\]
and Corollary \ref{c.int} implies that 
\begin{align}
\label{e.a}
F_n^\sigma(I_t^n(\alpha))
	&= \sum_{a=0}^{q_\alpha} b_\alpha^a t^a J_n(\tilde{f}_{\alpha,a}
		\hat{F}_n^{\sigma,\alpha})_t
\end{align}
is well-defined.  In particular, if $\alpha_m=1$, then 
$f_\alpha=f_\alpha(s)$ does not depend on $t$, and Corollary \ref{c.int} 
implies that $F_n^\sigma(I_t^n(\alpha))$ is a $\mathfrak{v}$-valued 
$L^2$-martingale.

The next two results show that $\hat{F}_n^{\sigma,\alpha}$ is
Hilbert-Schmidt as desired.

\begin{lem}
\label{l.HS}
Let $n\in\{2,\ldots,r\}$, $\sigma\in\mathcal{S}_n$, and
$\alpha\in\mathcal{J}_n$. 
For any $v\in\mathfrak{v}$, $\langle \hat{F}_n^{\sigma,\alpha}, v \rangle$ 
is a Hilbert-Schmidt operator on $\mathfrak{g}_{CM}^{\otimes p_\alpha}$.
\end{lem}

\begin{proof}
First consider the case $n=2$. In this case, $p_\alpha=0$ or $p_\alpha=2$.
If $p_\alpha=0$, then $\hat{F}_2^{\sigma,\alpha}=\sum_{i=1}^\infty
F_2^\sigma(h_i\otimes h_i)=0$.
If $p_\alpha=2$, then $\hat{F}_2^{\sigma,\alpha}(k_1\otimes
k_2)=F_2^\sigma(k_1\otimes k_2)=[k_{\sigma(1)},k_{\sigma(2)})]$ is
Hilbert-Schmidt by Corollary \ref{c.normest}, and thus $\langle
\hat{F}_2^{\sigma,\alpha}, v \rangle$ is Hilbert-Schmidt.
For $n=3$, $p_\alpha=1$ or $p_\alpha=3$.
If $p_\alpha=3$, then $\alpha=(1,1,1)$ and
\[ \hat{F}_3^{\sigma,\alpha}(k_1\otimes k_2\otimes k_3) 
	= F_3^{\sigma'}(k_1\otimes k_2\otimes k_3) 
	=[ [k_{\sigma(1)},k_{\sigma(2)}],k_{\sigma(3)}] \]
is Hilbert-Schmidt, again by Corollary \ref{c.normest}.  If $p_\alpha=1$, then
$\alpha=(1,2)$ or $\alpha=(2,1)$ and  
\[ \hat{F}_3^{\sigma,\alpha}(k) 
	= \sum_{i=1}^\infty F_3^{\sigma'}(k\otimes h_i\otimes h_i),
\] 
and we need only consider the case that
\[ F_3^{\sigma'}(k\otimes h\otimes h) = [[h,k],h]. \]
So let $\{k_i\}_{i=1}^\infty$ be an orthonormal
basis of $\mathfrak{g}_{CM}$ and $\{e_\ell\}_{\ell=1}^N$ be an orthonormal
basis of $\mathfrak{v}$.  As in the proof of Corollary \ref{c.normest}, 
expanding terms in an orthonormal basis of $\mathfrak{v}$
and applying the Cauchy-Schwarz inequality gives
\begin{align*}
\|\langle \hat{F}_3^{\sigma,\alpha}, v \rangle\|_2^2 
	&=\sum_{i=1}^\infty \left|\sum_{j=1}^\infty 
		\langle [[h_j,k_i],h_j],v\rangle\right|^2 
	= \sum_{i=1}^\infty \left|\sum_{j=1}^\infty 
		\sum_{\ell=1}^N \langle [e_\ell,h_j], v\rangle
		\langle e_\ell,[h_j,k_i]\rangle \right|^2 \\
	&\le N \sum_{i=1}^\infty 
		\sum_{\ell=1}^N \left|\sum_{j=1}^\infty \langle [e_\ell,h_j],v\rangle
		\langle e_\ell,[h_j,k_i]\rangle \right|^2 \\
	&\le N \sum_{i=1}^\infty \sum_{\ell=1}^N \left(\sum_{j=1}^\infty 
		|\langle [e_\ell,h_j],v\rangle|^2\right)
		\left(\sum_{j=1}^\infty |\langle e_\ell,[h_j,k_i]\rangle|^2\right) \\
	&\le N \left( \sum_{j=1}^\infty \sum_{\ell=1}^N 
		|\langle [e_\ell,h_j],v\rangle|^2\right)
		\left(\sum_{i,j=1}^\infty \sum_{\ell=1}^N 
		|\langle e_\ell,[h_j,k_i]\rangle|^2\right) \\
	&\le N\|v\|^2<D-1> \|[\cdot,\cdot]\|_2^2 \cdot \|[\cdot,\cdot]\|_2^2.
\end{align*}

Now assume $\langle \hat{F}_{n-1}^{\sigma,\alpha}, v \rangle$ is
Hilbert-Schmidt for all $\sigma\in\mathcal{S}_{n-1}$ and
$\alpha\in\mathcal{J}_{n-1}$, and consider $\langle
\hat{F}_n^{\sigma,\alpha}, v \rangle$ for some $\sigma\in\mathcal{S}_n$
and $\alpha\in\mathcal{J}_n^m$.  Let $a=p_\alpha$ and
$b=q_\alpha$, and note that either $a\ge1$ and
\begin{align}
\hat{F}_n^{\sigma,\alpha}&(k_1\otimes\cdots\otimes k_a) \notag\\
	&= \sum_{j_1,\ldots,j_b=1}^\infty 
		F_n^{\sigma'}(k_1\otimes\cdots\otimes k_a
		\otimes h_{j_1}\otimes h_{j_1}
		\otimes\cdots\otimes h_{j_b}\otimes h_{j_b}) \notag\\
	&= \sum_{j_1,\ldots,j_b=1}^\infty  
		[F_{n-1}^{\sigma''}(k_1\otimes\cdots\otimes k_{d-1}\otimes
		k_{d+1}\otimes\cdots\otimes k_a\otimes 
		h_{j_1}\otimes\cdots\otimes h_{j_b}),k_d] \notag\\
	&\label{e.case1}
	= [\hat{F}_{n-1}^{\tau,\beta}(k_1\otimes\cdots\otimes k_{d-1} \otimes
		k_{d+1}\otimes\cdots\otimes k_a), k_d],
\end{align}
for some $d\in\{1,\ldots,a\}$, $\sigma'',\tau\in\mathcal{S}_{n-1}$, 
and $\beta\in\mathcal{J}_{n-1}^{m-1}$ such that $p_\beta=p_\alpha-1$ and
$q_\beta=q_\alpha$, or $b\ge1$ and
\begin{align}
&\hat{F}_n^{\sigma,\alpha}(k_1\otimes\cdots\otimes k_a) \notag \\
	&\,= \sum_{j_1,\ldots,j_b=1}^\infty 
		[F_{n-1}^{\sigma''}(k_1\otimes\cdots\otimes k_a\otimes 
		h_{j_1}\otimes\cdots\otimes h_{j_{d-1}}
		\otimes h_{j_d}\otimes h_{j_{d+1}}
		\otimes\cdots\otimes h_{j_b}),h_{j_d}] \notag \\
	&\label{e.case2}
	\,= \sum_{j_d=1}^\infty
		[\hat{F}_{n-1}^{\tau,\beta}(k_1\otimes\cdots\otimes k_a\otimes 
		h_{j_d}), h_{j_d}],
\end{align}
for some $d\in\{1,\ldots,b\}$, $\sigma'',\tau \in\mathcal{S}_{n-1}$
and $\beta\in\mathcal{J}_{n-1}^m$ such that $p_\beta=p_\alpha+1$ and
$q_\beta=q_\alpha-1$.  In the first case, working as above for $n=3$,
\begin{align*}
&\left\|\langle \hat{F}_n^{\sigma,\alpha},v \rangle \right\|_2^2
	= \sum_{i_1,\ldots,i_a=1}^\infty \left|
		\sum_{j_1,\ldots,j_b=1}^\infty \langle F_n^{\sigma'}(k_{i_1}
		\otimes\cdots\otimes k_{i_a}\otimes 
		h_{j_1}\otimes\cdots\otimes h_{j_b}), v\rangle \right|^2  \\
	&\quad= \sum_{i_1,\ldots,i_a=1}^\infty \left|
		\sum_{j_1,\ldots,j_b=1}^\infty 
		\langle[F_{n-1}^{\sigma''}(k_{i_1}\otimes\cdots\otimes
		h_{j_b}),k_{i_d}],v\rangle \right|^2 \\
	&\quad\le N \sum_{i_1,\ldots,i_a=1}^\infty \sum_{\ell=1}^N
		\left| \sum_{j_1,\ldots,j_b=1}^\infty 		
		\langle F_{n-1}^{\sigma''}(k_{i_1}\otimes\cdots\otimes
		h_{j_b}),e_\ell\rangle \langle[e_\ell,k_{i_d}],v\rangle \right|^2 \\
	&\quad= N \sum_{i_1,\ldots,i_a=1}^\infty \sum_{\ell=1}^N
		|\langle[e_\ell,k_{i_d}],v\rangle|^2 
		\left| \sum_{j_1,\ldots,j_b=1}^\infty 		
		\langle F_{n-1}^{\sigma''}(k_{i_1}\otimes\cdots\otimes
		h_{j_b}),e_\ell\rangle \right|^2 \\
	&\quad\le N \|v\|^2 \|[\cdot,\cdot]\|^2_2
		\sum_{\ell=1}^N \left\|\langle \hat{F}_{n-1}^{\tau,\beta},
		e_\ell\rangle\right\|_2^2,
\end{align*}
which is finite by the induction hypothesis.  Similarly, in the second case
\begin{align*}
&\left\|\langle \hat{F}_n^{\sigma,\alpha},v \rangle \right\|_2^2
	= \sum_{i_1,\ldots,i_a=1}^\infty \left|
		\sum_{j_1,\ldots,j_b=1}^\infty 
		\langle[F_{n-1}^{\sigma''}(k_{i_1}\otimes\cdots\otimes
		h_{j_b}),h_{j_d}],v\rangle \right|^2 \\
	&\quad\le N \sum_{i_1,\ldots,i_a=1}^\infty \sum_{\ell=1}^N
		\left| \sum_{j_1,\ldots,j_b=1}^\infty 		
		\langle F_{n-1}^{\sigma''}(k_{i_1}\otimes\cdots\otimes
		h_{j_b}),e_\ell\rangle \langle[e_\ell,h_{j_d}],v\rangle \right|^2 \\
	&\quad\le N \left(\sum_{i_1,\ldots,i_a=1}^\infty \sum_{\ell=1}^N
		\sum_{j_d=1}^\infty \left| 
		\sum_{j_1,\ldots,j_{d-1},j_{d+1},\ldots,j_b=1}^\infty 		
		\langle F_{n-1}^{\sigma'}(k_{i_1}\otimes\cdots\otimes
		h_{j_b}),e_\ell\rangle\right|^2 \right) \\
	&\qquad\qquad \times\left(\sum_{\ell=1}^N \sum_{j_d=1}^\infty 
		|\langle[e_\ell,h_{j_d}],v\rangle|^2\right) \\
	&\quad\le N \sum_{\ell=1}^N \left\|\langle \hat{F}_{n-1}^{\tau,\beta},
		e_\ell\rangle\right\|_2^2 \cdot \|v\|^2\|[\cdot,\cdot]\|_2^2.
\end{align*}  
\end{proof}

\begin{prop}
\label{p.HS}
Let $n\in\{2,\ldots,r\}$, $\sigma\in\mathcal{S}_n$, and 
$\alpha\in\mathcal{J}_n$.
Then $\hat{F}_n^{\sigma,\alpha}:\mathfrak{g}_{CM}^{\otimes p_\alpha}
\rightarrow \mathfrak{v}$ is Hilbert-Schmidt.
\end{prop}
\begin{proof}
This proof is analogous to that of Lemma \ref{l.HS}.
For $\hat{F}_n^{\sigma,\alpha}$ as in equation (\ref{e.case2}), we have
\begin{align*}
\|\hat{F}_n^{\sigma,\alpha}&\|_2^2
	= \sum_{i_1,\ldots,i_a=1}^\infty \left\|
		\sum_{j_1,\ldots,j_b=1}^\infty 
		[F_{n-1}^{\sigma''}(k_{i_1}\otimes\cdots\otimes 
		h_{j_b}),h_{j_d}] \right\|^2 \\
	&\le N \sum_{i_1,\ldots,i_a=1}^\infty \sum_{\ell=1}^N
		\left\| \sum_{j_1,\ldots,j_b=1}^\infty 		
		\langle F_{n-1}^{\sigma''}(k_{i_1}\otimes\cdots\otimes 
		h_{j_b}), e_\ell\rangle [e_\ell,h_{j_d}] \right\|^2 \\
	&\le N \left(\sum_{\ell=1}^\infty \sum_{j_\ell=1}^\infty 
		\|[e_\ell,h_{j_d}]\|^2 \right) \\
	&\qquad \times
		\left(\sum_{i_1,\ldots,i_a=1}^\infty \sum_{\ell=1}^N
		\sum_{j_d=1}^\infty \left| 
		\sum_{j_1,\ldots,j_{d-1},j_{d+1},\ldots,j_b=1}^\infty 		
		\langle F_{n-1}^{\sigma''}(k_{i_1}\otimes\cdots\otimes
		h_{j_b}), e_\ell\rangle \right|^2\right) \\
	&\le N \|[\cdot,\cdot]\|_2^2 \sum_{d=1}^N 
		\left\|\langle \hat{F}_{n-1}^{\tau,\beta}
		,e_\ell\rangle \right\|_2^2,
\end{align*}
which is finite by Corollary \ref{c.normest} and Lemma \ref{l.HS}.  
In a similar way, one may show 
that, for $\hat{F}_n^{\sigma,\alpha}$ as in equation (\ref{e.case1}),
\[ \|\hat{F}_n^{\sigma,\alpha}\|_2^2
	\le N \|[\cdot,\cdot]\|_2^2 \sum_{d=1}^N
		\left\|\langle \hat{F}_{n-1}^{\tau,\beta},
		e_\ell\rangle \right\|_2^2. \]
\end{proof}

\begin{remark}
The proofs of the previous propositions rely strongly on $\mathfrak{v}$ being
finite dimensional.  Thus, if we were to extend the results of this paper to
$\mathfrak{v}$ an infinite dimensional Lie algebra, another
proof would be required here, or more likely, some trace class requirements
on the Lie bracket of $\mathfrak{g}$.
\end{remark}

Proposition \ref{p.HS} allows us to make the following definition.

\begin{defn}
\label{d.BM}
A {\em Brownian motion} on $G$ is the continuous $G$-valued process defined by 
\[g_t = \sum_{n=1}^r \sum_{\sigma\in\mathcal{S}_n}
		\sum_{m=\lceil n/2\rceil}^n \frac{c^\sigma_n }{2^{n-m}}
		\sum_{\alpha\in\mathcal{J}_n^m}
		\int_{\Delta_{p_\alpha}(t)} f_\alpha(s,t) \hat{F}_n^{\sigma,\alpha}
		(dB_{s_1}\otimes\cdots\otimes dB_{s_{p_\alpha}}), \]
where 
\[ c_n^\sigma=(-1)^{e(\sigma)}\bigg/n^2\begin{bmatrix} n-1 \\ e(\sigma) 
	\end{bmatrix}, \]
$\hat{F}_n^{\sigma,\alpha}$ is as defined in (\ref{e.hatF}) and $f_\alpha$ is
a polynomial of degree $q_\alpha$ in $s=(s_1,\ldots,s_{p_\alpha})$ and $t$ as
described in (\ref{e.falpha}).
For $t>0$, let $\nu_t=\mathrm{Law}(g_t)$ be the {\em heat kernel measure at
time $t$}, a probability measure on $G$.
\end{defn}

\begin{example}[The step 3 case]
Suppose that $\mathfrak{g}$ is nilpotent of step 3.  Then
\begin{align*} 
g_t &= \sum_{n=1}^3 \sum_{\sigma\in\mathcal{S}_n} c_n^\sigma
		F_n^\sigma(M_n(t)) \\
	&= \sum_{n=1}^3 \sum_{\sigma\in\mathcal{S}_n}
	\sum_{m=\lceil n/2\rceil}^n \frac{c^\sigma_n }{2^{n-m}}
		\sum_{\alpha\in\mathcal{J}_n^m} F_n^\sigma(I_t^n(\alpha))\\
	&= \sum_{n=1}^3 \sum_{\sigma\in\mathcal{S}_n}
		\sum_{m=\lceil n/2\rceil}^n \frac{c^\sigma_n }{2^{n-m}}
		\sum_{\alpha\in\mathcal{J}_n^m}
		\int_{\Delta_{p_\alpha}(t)} f_\alpha(s,t) \hat{F}_n^{\sigma,\alpha}
		(dB_{s_1}\otimes\cdots\otimes dB_{s_{p_\alpha}}).
\end{align*}
For $n=1$, there is the single term given by 
\[ M_1(t)=\int_0^t \delta B_s = B_t. \]
For $n=2$, $\mathcal{J}_2=\{(1,1),(2)\}$, and so
\begin{align*}
M_2(t) 
	&= I_t^2((1,1)) + \frac{1}{2}I_t^2((2)) \\
	&= \int_{\Delta_2(t)} dB_{s_1} \otimes dB_{s_2}
		+ \frac{1}{2} \int_0^t h_i \otimes h_i \,ds_2 \\
	&= \int_{\Delta_2(t)} dB_{s_1} \otimes dB_{s_2}
		+ \frac{1}{2}t \sum_{i=1}^\infty h_i \otimes h_i.
\end{align*}
There are of course just two permutations: $\sigma=(12)$
with $e(\sigma)=0$ and $c_2^\sigma=\frac{1}{4}$, and $\tau=(21)$ with
$e(\tau)=1$ and  $c_2^\tau=-\frac{1}{4}$, and, by the antisymmetry of the Lie
bracket,
\[ \sum_{\sigma\in\mathcal{S}_2} c_2^\sigma F_3^\sigma(M_2(t))
	=  \frac{1}{4}[dB_{s_1},dB_{s_2}] - \frac{1}{4}[dB_{s_2},dB_{s_1}] 
	= \frac{1}{2} [dB_{s_1},dB_{s_2}]. \]

For $n=3$, the permutations are $(123)$ with $e=0$, $(213)$, $(132)$, $(312)$,
$(231)$ with $e=1$, and $(321)$ with $e=2$.
Thus,
\begin{align}
\sum_{\sigma\in\mathcal{S}_3} c_3^\sigma F_3^\sigma(k_1\otimes k_2\otimes k_3)
	&= \frac{1}{9}    [[k_1,k_2],k_3] - \frac{1}{18} [[k_2,k_1,],k_3]
		- \frac{1}{18} [[k_1,k_3],k_2] \notag\\
	&\qquad - \frac{1}{18} [[k_3,k_1],k_2] - \frac{1}{18} [[k_2,k_3],k_1] 
		+ \frac{1}{9}  [[k_3,k_2,],k_1] \notag\\
	&\label{e.blah}
	= \frac{1}{6} [[k_1,k_2],k_3] + \frac{1}{6} [[k_3,k_2,],k_1].
\end{align}
Also, $\mathcal{J}_3=\{(1,1,1),(1,2),(2,1)\}$, and so
\begin{align*}
M_3(t) 
	&= I_t^3((1,1,1)) + \frac{1}{2}I_t^3((1,2)) + \frac{1}{2} I_t^3((2,1))\\
	&= \int_{\Delta_3(t)} dB_{s_1} \otimes dB_{s_2} \otimes dB_{s_3} 
		+ \frac{1}{2} \int_{\Delta_2(t)} \sum_{i=1}^\infty
		dB_{s_1} \otimes h_i\otimes h_i ds_3 \\
	&\qquad + \frac{1}{2} \int_0^t \sum_{i=1}^\infty
		s_3 h_i \otimes h_i \otimes dB_{s_3}  \\
	&= \int_{\Delta_3(t)} dB_{s_1} \otimes dB_{s_2} \otimes dB_{s_3} 
		+ \frac{1}{2} \int_0^t \sum_{i=1}^\infty
		(t-s_1) dB_{s_1}\otimes  h_i\otimes h_i  \\
	&\qquad + \frac{1}{2} \int_0^t \sum_{i=1}^\infty
		s_3 h_i \otimes h_i \otimes dB_{s_3}.  
\end{align*}
Note that $f_{(1,2)}(s,t)=t-s_1$ and $f_{(2,1)}(s,t)=s_3$.  Plugging this into
equation (\ref{e.blah}) gives, for the $\alpha=(1,1,1)\in\mathcal{J}_3^3$ term,
\begin{align*} 
\sum_{\sigma\in\mathcal{S}_3} c_3^\sigma F_3^\sigma(I_t^3((1,1,1)))
	&= \sum_{\sigma\in\mathcal{S}_3} c_3^\sigma \int_{\Delta_3(t)}
		F_3^\sigma (dB_{s_1}\otimes dB_{s_2}\otimes dB_{s_3}) \\
	&= \frac{1}{6} \int_{\Delta_3(t)}([[dB_{s_1},dB_{s_2}], dB_{s_3}]
		+ [[dB_{s_3},dB_{s_2}], dB_{s_1}] ). 
\end{align*}
For $\alpha=(1,2)\in\mathcal{J}_3^2$,
\[ \sum_{\sigma\in\mathcal{S}_3} c_3^\sigma F_3^\sigma(I_t(1,2))
	= \frac{1}{6} \int_0^t \sum_{i=1}^\infty (t-s_1)[[dB_{s_1},h_i], h_i], \]
and 
\[ \hat{F}_3^{\sigma,(1,2)}(k) = \sum_{i=1}^\infty F_3^\sigma(k\otimes
	h_i\otimes h_i ) \]
with $\sigma'=\sigma$. For $\alpha=(2,1)\in\mathcal{J}_3^2$,
\[ \sum_{\sigma\in\mathcal{S}_3} c_3^\sigma F_3^\sigma (I_t((2,1)))
	= \frac{1}{6} \int_0^t \sum_{i=1}^\infty s_3 [[dB_{s_3},h_i], h_i], \]
and note that, in this case,
\[ \hat{F}_3^{\sigma,(2,1)}(k) = \sum_{i=1}^\infty
		F_3^{\sigma'}(k\otimes h_i\otimes h_i ) 
	= \sum_{i=1}^\infty F_3^{\sigma}(h_i\otimes h_i\otimes k),\]
and so $\sigma'=\sigma\circ (231)$ (or $\sigma'=\sigma\circ(321)$).
Combining the above, Brownian motion on $G$ may be written as
\begin{align*}
g_t &= B_t + \frac{1}{2}\int_{\Delta_2(t)} [dB_{s_1},dB_{s_2}] \\
	&\quad + \frac{1}{12}\int_{\Delta_3(t)} ([[dB_{s_1},dB_{s_2}], dB_{s_3}]
		+ [[dB_{s_3},dB_{s_2}], dB_{s_1}]) \\
	&\quad + \frac{1}{24}\sum_{i=1}^\infty\int_0^t ( (t-s) [[dB_s,h_i], h_i] 
		+ s [[dB_s,h_i], h_i] ) \\
	&= B_t + \frac{1}{2}\int_0^t [B_s,dB_s] 
		+ \frac{1}{12}\int_{\Delta_2(t)} ([[B_{s_1},dB_{s_1}], dB_{s_2}]
		+ [[dB_{s_2},dB_{s_1}], B_{s_1}]) \\
	&\quad + \frac{1}{24}\sum_{i=1}^\infty t[[B_t,h_i], h_i].
\end{align*}
\end{example}

\begin{remark}
In principle, the Brownian motion on $G$ has generator 
\[ \Delta = \sum_{i=1}^\infty \tilde{h}_i^2, \]
where $\{h_i\}_{i=1}^\infty$ is an orthonormal basis of
$\mathfrak{g}_{CM}=H\oplus\mathfrak{v}$ and $\tilde{h}$ is the
unique left invariant vector field on $G$ such that
$\tilde{h}(\mathbf{e})=h$, and $\Delta$ is well-defined independent
of the choice of orthonormal basis.  Then the heat kernel measure
$\{\nu_t\}_{t>0}$ has the standard characterization as the unique family of
probability measures such that $\nu_t(f):=\int_G f \,d\nu_t$
is continuously differentiable in $t$ for all $f\in
C_b^2(G)$ and satisfies
\[ \frac{d}{dt} \nu_t(f) = \frac{1}{2}\nu_t(\Delta f) \quad \text{ with }
	\lim_{t\downarrow0}\nu_t(f)=f(e). \]
However, this realization of $\nu_t$ is not necessary for our results.
\end{remark}

\begin{prop}[Finite dimensional approximations]
\label{p.approx}
For $P\in\mathrm{Proj}(W)$, let $g^P_t$ be the continuous process on $G_P$
defined by
\[ g_t^P = \sum_{n=1}^r \sum_{\sigma\in\mathcal{S}_n}
	\sum_{m=\lceil n/2\rceil}^n \frac{c^\sigma_n }{2^{n-m}}
		\sum_{\alpha\in\mathcal{J}_n^m}
	\int_{\Delta_{p_\alpha}(t)} f_\alpha(s,t) \hat{F}_n^{\sigma,\alpha}
	(d\pi B_{s_1}\otimes\cdots\otimes d\pi B_{s_{p_\alpha}}), \]
for $\pi(w,x)=(Pw,x)$.
Then $g_t^P$ is Brownian motion on $G_P$.  In particular, let
$g^\ell_t=g_t^{P_\ell}$, for projections
$\{P_\ell\}_{\ell=1}^\infty\subset\mathrm{Proj}(W)$ 
as in equation (\ref{e.pm}).  Then, for all $p\in[1,\infty)$ and $t<\infty$,
\begin{equation} 
\label{e.b}
\lim_{\ell\rightarrow\infty}\mathbb{E}\left[\sup_{\tau\le t}
	\left\|g^\ell_\tau-g_\tau\right\|_\mathfrak{g}^p\right] = 0. 
\end{equation}
\end{prop}
\begin{proof}
First note that $g_t^P$ solves the Stratonovich equation 
$\delta g_t^P = L_{g_t^P*}\delta PB_t$ with $g_0^P=\mathbf{e}$, see
\cite{BenArous89,Castell93,Baudoin04}.  Thus, $g_t^P$ is a $G_P$-valued
Brownian motion.  

Now, if $\beta_t$ a Brownian motion on $W$, then, for all $p\in[1,\infty)$,
\[ \lim_{\ell\rightarrow\infty}\mathbb{E}\left[\sup_{\tau\le t}
	\left\|P_\ell \beta_\tau - \beta_\tau\right\|_W^p\right] = 0; \]
see, for example, Proposition 4.6 of \cite{DG08-2}.  Thus,
\[ \lim_{\ell\rightarrow\infty}\mathbb{E}\left[\sup_{\tau\le t}
	\left\|\pi_\ell B_\tau - B_\tau\right\|_\mathfrak{g}^p\right] = 0. \]

By equation (\ref{e.a}) and its preceding discussion,
\[ g_t^\ell = \sum_{n=1}^r \sum_{\sigma\in\mathcal{S}_n}
	\sum_{m=\lceil n/2\rceil}^n \frac{c^\sigma_n }{2^{n-m}}
		\sum_{\alpha\in\mathcal{J}_n^m}
	\sum_{a=0}^{q_\alpha} b_\alpha^a t^a J_n^\ell(\tilde{f}_\alpha
		\hat{F}_n^{\sigma,\alpha})_t, \]
and thus, to verify (\ref{e.b}), it suffices to show that, for all
$p\in[1,\infty)$,
\[ \lim_{\ell\rightarrow\infty} \mathbb{E}\left[ \sup_{\tau\le t} 
	\left\|J_n^\ell(\tilde{f}_\alpha\hat{F}_n^{\sigma,\alpha})_\tau -
	J_n(\tilde{f}_\alpha\hat{F}_n^{\sigma,\alpha})_\tau\right\|_\mathfrak{v}^p
	\right] = 0, \]
for all $n\in\{2,\ldots,r\}$, $\sigma\in\mathcal{S}_n$ and
$\alpha\in\mathcal{J}_n$.  By Proposition \ref{p.HS},
$\hat{F}_n^{\sigma,\alpha}$ is Hilbert-Schmidt, and recall that
$\tilde{f}_\alpha$ is a deterministic polynomial function in $s$.  Thus
$J_n^\ell(\tilde{f}_\alpha \hat{F}_n^{\sigma,\alpha})$ and
$J_n(\tilde{f}_\alpha \hat{F}_n^{\sigma,\alpha})$ are 
$\mathfrak{v}$-valued martingales as defined in Corollary \ref{c.int}.
So, by Doob's maximal inequality, it suffices to show that 
\[ \lim_{\ell\rightarrow\infty} \mathbb{E}
	\left\|J_n^\ell(\tilde{f}_\alpha\hat{F}_n^{\sigma,\alpha})_t -
	J_n(\tilde{f}_\alpha\hat{F}_n^{\sigma,\alpha})_t 
	\right\|_\mathfrak{v}^p  = 0  \] 
Corollary \ref{c.int} gives the limit for $p=2$.  For $p>2$, since each
$J_n^\ell(\tilde{f}_\alpha \hat{F}_n^{\sigma,\alpha})$ and
$J_n(\tilde{f}_\alpha \hat{F}_n^{\sigma,\alpha})$ has chaos expansion
terminating at degree $n$, a theorem of Nelson (see Lemma 2 of
\cite{Nelson73b} and pp. 216-217 of \cite{Nelson73c}) implies that,
for each $j\in\mathbb{N}$, there exists $c_j<\infty$ such that 
\[ \mathbb{E}\left\|J_n^\ell(\tilde{f}_\alpha\hat{F}_n^{\sigma,\alpha})_t -
		J_n(\tilde{f}_\alpha\hat{F}_n^{\sigma,\alpha})_t 
		\right\|_\mathfrak{v}^{2j}
	\le c_j \left(\mathbb{E}\left\|
		J_n^\ell(\tilde{f}_\alpha\hat{F}_n^{\sigma,\alpha})_t -
		J_n(\tilde{f}_\alpha\hat{F}_n^{\sigma,\alpha})_t 
		\right\|_\mathfrak{v}^2\right)^j. \]
\end{proof}

\section{Heat kernel measure}
\label{s.hki}

We collect here some properties of the heat kernel measure on $G$.
The following results are completely analogous to Corollary 4.9
of \cite{DG08-2} and Proposition 4.6 in \cite{DG08-3}.  
The proofs are included here for the convenience of the
reader.

\begin{prop}
\label{c.4.9}
For any $t>0$, the heat kernel measure $\nu_t$ is invariant under the
inversion map $g\mapsto g^{-1}$ for any $g\in G$.
\end{prop}
\begin{proof}
The heat kernel measures $\nu_t^{P_n}=\mathrm{Law}(g_t^n)$ on the finite dimensional groups
$G_{P_n}$ are invariant under inversion (see, for example,
\cite{Driver97}).  Suppose that $f:G\rightarrow\mathbb{R}$ is a
bounded continuous function.  By passing to a subsequence if
necessary, we may assume that the sequence of $G_{P_n}$-valued
random variables $\{g^n_t\}_{n=1}^\infty$ in Proposition \ref{p.approx}
converges almost surely to $g_t$.  Thus, by dominated convergence,
\[ \mathbb{E}\left[f\left(g_t^{-1}\right)\right]
	= \lim_{n\rightarrow\infty} 
		\mathbb{E}\left[f\left( (g^n_t)^{-1}\right)\right]
	= \lim_{n\rightarrow\infty} 
		\mathbb{E}\left[f\left(g^n_t\right)\right]
	= \mathbb{E}\left[f\left(g_t\right)\right].
\]
Since $\nu_t$ is the law of $g_t$, this completes the proof.
\end{proof}  

\begin{prop}
For all $t>0$, $\nu_t(G_{CM})=0$.
\end{prop}

\begin{proof}
Let $\mu_t$ denote Wiener measure on $W$ with variance $t$.
Then for a bounded measurable function $f$ on $G=W\oplus\mathfrak{v}$ such that
$f(w,x)=f(w)$,
\[ \int_G f(w)\,d\nu_t(w,x) = \mathbb{E}[f(\beta_t)]
	= \int_W f(w)\,d\mu_t(w). \]
Let $\pi:W\times\mathfrak{v}\rightarrow W$ be the projection $\pi(w,x)=w$.
Then $\pi_*\nu_t =\mu_t$, and thus
\[ \nu_t(G_{CM}) = \nu_t\left(\pi^{-1}(H)\right) = \pi_*\nu_t(H)
	= \mu_t(H) = 0. \] 
\end{proof}

This proposition gives some justification to our calling $G_{CM}$ the 
Cameron-Martin subgroup of $G$.  In the next section, we further justify this
by showing that a Cameron-Martin type quasi-invariance theorem holds 
for $\nu_t$.

\subsection{Quasi-invariance and Radon-Nikodym derivative estimates}
\label{s.quasi}

The following theorem states that the heat kernel measure $\nu_t =
\mathrm{Law}(g_t)$ is
quasi-invariant under left and right translation by elements of $G_{CM}$ and
gives estimates for the Radon-Nikodym derivatives of the translated measures.

\begin{thm}
\label{t.quasi}
For all $h\in G_{CM}$ and $t>0$, $\nu_t\circ L_h^{-1}$ and $\nu_t\circ
R_h^{-1}$ are absolutely continuous with respect to $\nu_t$.  Let 
\[ Z_h^l := \frac{d(\nu_t\circ L_h^{-1})}{d\nu_t} \qquad \text{ and }
	\qquad Z_h^r := \frac{d(\nu_t\circ R_h^{-1})}{d\nu_t} \]
be the Radon-Nikodym derivatives, $K$ be lower bound on the Ricci curvature of
$G$ as in Corollary \ref{c.Ric3}, and 
\[ c(t) := \frac{t}{e^t-1}, \qquad \text{ for all } t\in\mathbb{R}, \]
with the convention that $c(0)=1$.  Then, 
$Z_h^l,Z_h^r\in L^p(\nu_t)$ for all $p\in[1,\infty)$, 
and both satisfy the estimate
\[ \|Z_h^*\|_{L^p(\nu_t)} \le \exp \left( \frac{c(Kt)(p-1)}{2t}
	d^2_{CM}(\mathbf{e},h)\right), \]
where $*=l$ or $*=r$.
\end{thm}
\begin{proof}
As in \cite{DG08-2}, the proof of this theorem is an application of Theorem
7.3 and Corollary 7.4 in \cite{DG08-1} on the quasi-invariance of heat kernel
measures for inductive limits of finite dimensional Lie groups.  In applying
these results, the reader should take $G_0=G_{CM}$, $A=\mathrm{Proj}(W)$,
$s_P=\pi_P$, $\nu_P=\mathrm{Law}(g_t^P)$, and $\nu=\nu_t=\mathrm{Law}(g_t)$.
We now verify that the hypotheses of Theorem 7.3 of \cite{DG08-1} are
satisfied.  

By Corollary \ref{c.length}, the inductive
limit group $\cup_{P\in\mathrm{Proj}(W)} G_P$ is a
dense subgroup of $G_{CM}$.  By Proposition \ref{p.approx}, for any
$\{P_n\}_{n=1}^\infty\subset\mathrm{Proj}(W)$ with $P_n|_H\uparrow I_H$ and
$f:G\rightarrow\mathbb{R}$ a bounded continuous function,
\begin{equation}
\label{e.pf}
\int_G f\,d\nu = \lim_{n\rightarrow\infty} \int_{G_{P_n}} (f\circ
	i_{P_n})\,d\nu_{P_n},
\end{equation}
and thus the heat kernel measure is consistent on finite dimensional
projections of $G_{CM}$.
Corollary \ref{c.Ric3} says that $K>-\infty$ and
$\mathrm{Ric}^P\ge Kg^P$, for all
$P\in\mathrm{Proj}(W)$, and thus the Ricci curvature is uniformly bounded on
these projections.  Lastly, the length of a path in the inductive limit group
can be approximated by the lengths of paths in the finite dimensional
projections.  That is, for any
$P_0\in\mathrm{Proj}(W)$ and $\varphi\in C^1([0,1],G_{CM})$ with
$\varphi(0)=\mathbf{e}$, there exists an increasing sequence
$\{P_n\}_{n=1}^\infty\subset\mathrm{Proj}(W)$ such that $P_0\subset P_n$,
$P_n|_H\uparrow I_H$, and 
\[ \ell_{CM}(\varphi) = \lim_{n\rightarrow\infty}
	\ell_{G_{P_n}}(\pi_n\circ\varphi). \]
To see this, let $\varphi(t)=(A(t),a(t))$ be a path in $G_{CM}$, and
recall that, by equation (\ref{e.lg}),
\begin{align*}
\ell_{G_{P_n}}(\pi_n\circ\varphi)
	&=\int_0^1 \left\|\pi_n\varphi'(s) + \sum_{\ell=1}^{r-1}
		d_\ell \mathrm{ad}_{\pi_n\varphi(s)}^\ell \pi_n\varphi'(s)
		\right\|_{\mathfrak{g}_{CM}}\,ds \\
	&= \int_0^1 \sqrt{ \left\|P_nA'(s)\right\|^2_H 
		+ \left\| a'(s) + \sum_{\ell=1}^{r-1}
		d_\ell \mathrm{ad}_{\pi_n\varphi(s)}^\ell \pi_n\varphi'(s)
		\right\|_\mathfrak{v}^2}\,ds 
\end{align*}
Applying dominated convergence to this equation shows that (\ref{e.pf}) holds 
for any such choice of $P_n|_H\uparrow I_H$ such that $P_0\subset P_n$.
\end{proof}

We also have the usual strong converse to quasi-invariance of
$\nu_t$ under translations by elements in $G_{CM}$.

\begin{prop}
For $h\in G\setminus G_{CM}$ and $t>0$, $(\nu_t\circ
L_h^{-1})$ and $\nu_t$ are singular and $(\nu_t\circ
R_h^{-1})$ and $\nu_t$ are singular.
\end{prop}

\begin{proof}
Again, let $\mu_t$ denote Wiener measure on $W$ with variance $t$.
Let $h=(A,a)\in G\setminus G_{CM}$ with $A\in W\setminus H$ and 
$a\in\mathfrak{v}$.  Given a measurable subset $U\subset W$,
\[ \nu_t(U\times\mathfrak{v}) = P(\beta_t\in U) = \mu_t(U). \]
If $A\in W\setminus H$, $\mu_t(\cdot-A)$ and $\mu_t$
are singular; for example, see Corollary 2.5.3 of \cite{Bog98}.  Thus, there
are disjoint subsets $W_0$ and $W_1$ of $W$ such that
$\mu_t(W_0)=1=\mu_t(W_1-A)$.  Note that
\[ L_k^{-1}(U\times\mathfrak{v})=R_k^{-1}(U\times\mathfrak{v})
	= (U-A)\times\mathfrak{v}. \]
Thus, for $G_i:=W_i\times\mathfrak{v}$ for $i=0,1$, $G$ is the disjoint union
of $G_0$ and $G_1$, and $\nu_t(G_0)=\mu_t(W_0)=1$ while
\[ \nu_t\left(R_k^{-1}(G_1)\right) = \nu_t\left(L_k^{-1}(G_1)\right)
	= \nu_t( (W_1-A)\times\mathfrak{v}) = \mu_t(W_1-A) = 1. \]
\end{proof}

\begin{prop}
For all $h\in G_{CM}$ and $t>0$, $Z_h^r(g) = Z_{h^{-1}}^l(g^{-1})$.
\end{prop}

\begin{proof}
By Proposition \ref{c.4.9}, $\nu_t$ is invariant under inversions.  Thus
\begin{align*}
\int_G f(g\cdot h)\,d\nu_t(g) 
	&= \int_G f\left(g^{-1}\cdot h\right)\,d\nu_t(g) 
	= \int_G f\left(\left(h^{-1}\cdot g\right)^{-1}\right)\,d\nu_t(g) \\
	&= \int_G f\left(g^{-1}\right) Z_{h^{-1}}^l (g) \,d\nu_t(g) 
	= \int_G f(g) Z_{h^{-1}}^l \left(g^{-1}\right)\,d\nu_t(g).
\end{align*}
\end{proof}

\subsection{Logarithmic Sobolev inequality}
\label{s.logsob}

\begin{defn}
\label{d.cyl} 
A function $f:G\rightarrow\mathbb{R}$ is said to be a
{\it (smooth) cylinder function} if $f=F\circ\pi_P$ for some
$P\in\mathrm{Proj}(W)$ and
some (smooth) function $F:G_P\rightarrow\mathbb{R}$.  Also, $f$ is a 
{\it cylinder polynomial} if $f=F\circ\pi_P$ 
for $F$ a polynomial function on $G_P$.
\end{defn}

\begin{thm}
\label{t.logsob}
Given a cylinder polynomial $f$ on $G$, let
$\nabla f:G\rightarrow\mathfrak{g}_{CM}$ be the gradient of $f$,
the unique element of $\mathfrak{g}_{CM}$ such that
\[ \langle\nabla f(g), h \rangle_{\mathfrak{g}_{CM}} = \tilde{h}f(g)
	:= f'(g)(L_{g*}h_\mathbf{e}), \]
for all $h\in\mathfrak{g}_{CM}$.  
Then for $K$ as in Corollary \ref{c.Ric3},
\[ \int_G (f^2\ln f^2)\,d\nu_t -
		\left(\int_G f^2\,d\nu_t \right)\cdot\ln\left(\int_G
		f^2\,d\nu_t\right)
	\le 2\frac{1-e^{-Kt}}{K} \int_G \|\nabla f\|_{\mathfrak{g}_{CM}}^2
		\,d\nu_t. \]
\end{thm}

\begin{proof}
Following the method of Bakry and Ledoux applied to $G_P$ (see Theorem 2.9 of 
\cite{Driver96} for the case needed here) shows that
\[ \mathbb{E}\left[\left(f^2\ln f^2\right)\left(g^P_t\right)\right] 
		- \mathbb{E}\left[f^2\left(g^P_t\right)\right] 
		\ln\mathbb{E}\left[f^2\left(g_t^P\right)\right]
	\le 2 \frac{1 - e^{-K_pt}}{K_P} \mathbb{E}\left\|(\nabla^P
		f)\left(g^P_t\right)\right\|^2_{\mathfrak{g}_{CM}^P},
\]
for $K_P$ as in equation (\ref{e.pah}).  Since the function $x\mapsto
(1-e^{-x})/x$ is decreasing and $K\le K_P$ for all $P\in\mathrm{Proj}(W)$, 
this estimate also holds with $K_P$ replaced with $K$.  Now applying
Proposition \ref{p.approx} to pass to the limit as $P\uparrow I$ gives 
the desired result.
\end{proof}

\begin{remark}
It is desirable to state Theorem \ref{t.logsob} for a larger class of
functions in $L^2(\nu_t)$.  To do this, one must prove that the gradient operator $\nabla:L^2(\nu_t)\rightarrow
L^2(\nu_t)\otimes\mathfrak{g}_{CM}$ is closable. 
Unfortunately, Theorem \ref{t.quasi}
doesn't give good information on the dependence of the Radon-Nikodym
derivatives $Z_h^l$ and $Z_h^r$ on $h$, and so at this point we can't prove
the necessary integration by parts formulae to show closability.
\end{remark}

\bibliographystyle{amsplain}
\bibliography{biblio}

\providecommand{\bysame}{\leavevmode\hbox to3em{\hrulefill}\thinspace}
\providecommand{\MR}{\relax\ifhmode\unskip\space\fi MR }
\providecommand{\MRhref}[2]{%
  \href{http://www.ams.org/mathscinet-getitem?mr=#1}{#2}
}
\providecommand{\href}[2]{#2}
\begin{thebibliography}{10}

\bibitem{Malliavin06}
H{\'e}l{\`e}ne Airault and Paul Malliavin, \emph{Quasi-invariance of {B}rownian
  measures on the group of circle homeomorphisms and infinite-dimensional
  {R}iemannian geometry}, J. Funct. Anal. \textbf{241} (2006), no.~1, 99--142.
  \MR{MR2264248 (2008b:60119)}

\bibitem{AMR00}
Dmitri Alekseevsky, Peter~W. Michor, and Wolfgang~A.F. Ruppert,
  \emph{Extensions of {L}ie algebras}, 2000.

\bibitem{Baudoin04}
Fabrice Baudoin, \emph{An introduction to the geometry of stochastic flows},
  Imperial College Press, London, 2004. \MR{MR2154760 (2006f:60003)}

\bibitem{BenArous89}
G{\'e}rard Ben~Arous, \emph{Flots et s\'eries de {T}aylor stochastiques},
  Probab. Theory Related Fields \textbf{81} (1989), no.~1, 29--77. \MR{MR981567
  (90a:60106)}

\bibitem{Besse87}
Arthur~L. Besse, \emph{Einstein manifolds}, Ergebnisse der Mathematik und ihrer
  Grenzgebiete (3) [Results in Mathematics and Related Areas (3)], vol.~10,
  Springer-Verlag, Berlin, 1987. \MR{MR867684 (88f:53087)}

\bibitem{Bog98}
Vladimir~I. Bogachev, \emph{Gaussian measures}, Mathematical Surveys and
  Monographs, vol.~62, American Mathematical Society, Providence, RI, 1998.
  \MR{MR1642391 (2000a:60004)}

\bibitem{CarEil99}
Henri Cartan and Samuel Eilenberg, \emph{Homological algebra}, Princeton
  Landmarks in Mathematics, Princeton University Press, Princeton, NJ, 1999,
  With an appendix by David A. Buchsbaum, Reprint of the 1956 original.
  \MR{MR1731415 (2000h:18022)}

\bibitem{Castell93}
Fabienne Castell, \emph{Asymptotic expansion of stochastic flows}, Probab.
  Theory Related Fields \textbf{96} (1993), no.~2, 225--239. \MR{MR1227033
  (94g:60110)}

\bibitem{CorGrn90}
Lawrence~J. Corwin and Frederick~P. Greenleaf, \emph{Representations of
  nilpotent {L}ie groups and their applications. {P}art {I}}, Cambridge Studies
  in Advanced Mathematics, vol.~18, Cambridge University Press, Cambridge,
  1990, Basic theory and examples. \MR{MR1070979 (92b:22007)}

\bibitem{DG08-2}
B.~Driver and M.~Gordina, \emph{Heat kernel analysis on infinite-dimensional
  {H}eisenberg groups}, J. Funct. Anal. \textbf{255} (2008), no.~2, 2395--2461.

\bibitem{DG08-1}
\bysame, \emph{Integrated {H}arnack inequalities on {L}ie groups}, Submitted,
  arXiv:0711.4392 (2008).

\bibitem{DG08-3}
\bysame, \emph{Square integrable holomorphic functions on infinite-dimensional
  {H}eisenberg type groups}, to appear in Probab. Theory Related Fields (2008).

\bibitem{Driver97}
Bruce~K. Driver, \emph{Integration by parts and quasi-invariance for heat
  kernel measures on loop groups}, J. Funct. Anal. \textbf{149} (1997), no.~2,
  470--547. \MR{MR1472366 (99a:60054a)}

\bibitem{Driver96}
Bruce~K. Driver and Terry Lohrenz, \emph{Logarithmic {S}obolev inequalities for
  pinned loop groups}, J. Funct. Anal. \textbf{140} (1996), no.~2, 381--448.
  \MR{MR1409043 (97h:58176)}

\bibitem{duiskolk}
J.~J. Duistermaat and J.~A.~C. Kolk, \emph{Lie groups}, Universitext,
  Springer-Verlag, Berlin, 2000. \MR{MR1738431 (2001j:22008)}

\bibitem{Fuks}
D.~B. Fuks, \emph{Cohomology of infinite-dimensional {L}ie algebras},
  Contemporary Soviet Mathematics, Consultants Bureau, New York, 1986,
  Translated from the Russian by A. B. Sosinski{\u\i}. \MR{MR874337
  (88b:17001)}

\bibitem{Inahama04}
Yuzuru Inahama, \emph{Logarithmic {S}obolev inequality for {$H\sp s\sb
  0$}-metric on pinned loop groups}, Infin. Dimens. Anal. Quantum Probab.
  Relat. Top. \textbf{7} (2004), no.~1, 1--26. \MR{MR2021645 (2005d:58061)}

\bibitem{Ito51}
Kiyosi It{\^o}, \emph{Multiple {W}iener integral}, J. Math. Soc. Japan
  \textbf{3} (1951), 157--169. \MR{MR0044064 (13,364a)}

\bibitem{Kuo75}
Hui~Hsiung Kuo, \emph{Gaussian measures in {B}anach spaces}, Lecture Notes in
  Mathematics, Vol. 463, Springer-Verlag, Berlin, 1975. \MR{MR0461643 (57
  \#1628)}

\bibitem{Nelson73c}
Edward Nelson, \emph{The free {M}arkoff field}, J. Functional Analysis
  \textbf{12} (1973), 211--227. \MR{MR0343816 (49 \#8556)}

\bibitem{Nelson73b}
\bysame, \emph{Quantum fields and {M}arkoff fields}, Partial differential
  equations ({P}roc. {S}ympos. {P}ure {M}ath., {V}ol. {XXIII}, {U}niv.
  {C}alifornia, {B}erkeley, {C}alif., 1971), Amer. Math. Soc., Providence,
  R.I., 1973, pp.~413--420. \MR{MR0337206 (49 \#1978)}

\bibitem{Sadetov05}
S.~T. Sad{\`e}tov, \emph{On extensions of {L}ie algebras by the {H}eisenberg
  algebra}, Mat. Zametki \textbf{78} (2005), no.~5, 745--747. \MR{MR2252954
  (2007e:17016)}

\bibitem{Shigekawa04}
Ichiro Shigekawa, \emph{Stochastic analysis}, Translations of Mathematical
  Monographs, vol. 224, American Mathematical Society, Providence, RI, 2004,
  Translated from the 1998 Japanese original by the author, Iwanami Series in
  Modern Mathematics. \MR{MR2060917 (2005k:60002)}

\bibitem{Strichartz87}
Robert~S. Strichartz, \emph{The {C}ampbell-{B}aker-{H}ausdorff-{D}ynkin formula
  and solutions of differential equations}, J. Funct. Anal. \textbf{72} (1987),
  no.~2, 320--345. \MR{MR886816 (89b:22011)}

\bibitem{Varadarajan}
V.~S. Varadarajan, \emph{Lie groups, {L}ie algebras, and their
  representations}, Graduate Texts in Mathematics, vol. 102, Springer-Verlag,
  New York, 1984, Reprint of the 1974 edition. \MR{85e:22001}

\bibitem{Yankosky01}
Bill Yankosky, \emph{On nilpotent extensions of {L}ie algebras}, Houston J.
  Math. \textbf{27} (2001), no.~4, 719--724. \MR{MR1874666 (2002k:17026)}

\end{thebibliography}

\end{document}